\documentclass[12pt,fleqn,leqno]{article}
\usepackage{graphicx,color,enumerate,amssymb}
\usepackage{amsmath}
\usepackage{amssymb}
\usepackage{epsfig}

\makeatletter
\@addtoreset{equation}{section}
\renewcommand\thetable{\@arabic\c@table}
\renewcommand\thefigure{\@arabic\c@figure}
\long\def\@makecaption#1#2{%
  \vskip\abovecaptionskip
  \begin{center}%
  \sbox\@tempboxa{#1: #2}%
  \ifdim \wd\@tempboxa >\hsize
    #1: #2\par
  \else
    \global \@minipagefalse
    \hb@xt@\hsize{\hfil\box\@tempboxa\hfil}%
  \fi
  \end{center}%
  \vskip\belowcaptionskip}
\def\N{{\rm I\kern-.15em N}}
\def\R{{\rm I\kern-.2em R}}
\def\Z{{\rm Z\kern-.26em Z}}
\makeatother
\newtheorem{thm}{Theorem}[section]

\newtheorem{rem}[thm]{Remark}
\newtheorem{cor}[thm]{Corollary}
\newtheorem{prop}[thm]{Proposition}

\newcommand{\be}{\begin{eqnarray}}
\newcommand{\ee}{\end{eqnarray}}
\newcommand{\bq}{\begin{eqnarray*}}
\newcommand{\eq}{\end{eqnarray*}}

\newcommand{\ii}{\textrm{i}}
\def \argmax{\mathop{\hbox{\rm arg max}}}

\newcommand{\bewend}{\hspace*{2mm}\rule{3mm}{3mm}}

\newcommand{\RR}{\mathbb{R}}
\newcommand{\PP}{\mathbb{P}}
\newcommand{\NN}{\mathbb{N}}
\newcommand{\BE}{\mathbb{E}}
\newcommand{\e}{\textrm{e}}
\newcommand{\LL}{{\textrm L}^2}
\newcommand{\vertk}{\stackrel{{\cal D}}{\longrightarrow}}
\newcommand{\stk}{\stackrel{{\mathbb{P}}}{\longrightarrow}}

\usepackage{color}

\begin{document}
\begin{center}
{ \LARGE\sc Characterizations of multinormality and corresponding tests of fit, including for garch models \\
} \vspace*{0.5cm} {\large\sc Norbert Henze}$^{1}$, {\large\sc M.D. Jim\'enez--Gamero}$^{2}$, {\large\sc Simos G. Meintanis}$^{3,4}$\\ \vspace*{0.5cm}
{\it $^{1}$Institute of Stochastics, Karlsruhe Institute of Technology, Karlsruhe, Germany \\
{\it $^{2}$Department of Statistics and Operations Research, University of Seville, Seville, Spain}\\
{\it $^{3}$Department of Economics, National and Kapodistrian University
of Athens, Athens, Greece} \\
{\it{$^{\rm{4}}$Unit for Business Mathematics and Informatics, North--West University, Potchefstroom, South Africa}}}\footnote{On sabbatical leave from the University of Athens}
\\
\end{center}

{\small {\bf Abstract.} We provide novel characterizations of multivariate normality that incorporate both the characteristic function and the moment generating function, and we employ these results to construct a class of affine invariant, consistent and easy-to-use goodness-of-fit tests for normality.
The test statistics are suitably weighted $L^2$-statistics, and we provide their asymptotic behavior both for i.i.d. observations as well as in the context of testing that the innovation distribution of a multivariate GARCH model is Gaussian. We also study the finite-sample behavior of the new tests and compare the new criteria with alternative existing tests.}
\vspace*{0.3cm} {AMS
2000 classification numbers:} 62H15, 62G20.
\\ \\


\section{Introduction}\label{sec_intro}
Let $X$ be a (univariate) random variable with characteristic function (CF) $\varphi_X(t) = \BE[\exp(\textrm{i}tX)]$, $t \in \RR$. Assuming that the moment generating function (MGF) $M_X(t) = \BE[\exp(tX)]$ of $X$ exists for each $t\in \RR$,   \cite{vo:14} proved that the identity
\begin{equation}\label{volk} \varphi_X(t)M_X(t)-1=0 \quad \textrm{for each } t\in \RR
\end{equation}
characterizes the zero-mean Gaussian distribution. Notice that, by assuming the existence of the MGF in an interval around zero, the moments of any order are implicitly supposed to exist. A further point is that a CF satisfying (\ref{volk}) is necessarily real-valued. Hence, the symmetry  of the distribution of $X$ around zero is also implicit in (\ref{volk}), see \cite{lu:70}, \S3.1.

In the following we extend the characterization in (\ref{volk}) to non-centered random variables and to random vectors of arbitrary dimension.
 We then construct a goodness-of-fit test for multivariate normality based on a suitable weighted $L^2$-statistic
     in which both the CF and the MGF are estimated nonparametrically. Several variants of this test criterion are also suggested.
     Furthermore, we work out the asymptotics of our test statistics  in the context of
  independent and identically distributed (i.i.d.) random vectors, thus obtaining a new test for multivariate normality, but also in the context
       of GARCH-type dependence. The latter test provides a novel method for assessing the celebrated question whether
       a Gaussian GARCH driven volatility process is adequate for explaining the heavy tails that are often observed with financial time series.

       The rest of the paper unfolds as follows. In Section \ref{sec_charac} we prove the above-mentioned extension of (\ref{volk}).
       On the basis of the resulting characterization, Section \ref{sec_teststat} suggests two new classes of affine invariant and easily computable statistics for testing for multivariate normality.
       Section \ref{secgaminfty} shows that a 'certain limit statistic' of one of these classes is a linear combination of two well-known measures of multivariate skewness, while the other is related to Mardia's multivariate sample kurtosis.
       In Section \ref{sec_iid} we derive the limit null distribution of the new test statistics in the i.i.d. setting. Section \ref{sec_consist} addresses the question of consistency of the new tests against general alternatives. Section \ref{sec_garch} considers this criterion in the context of multivariate GARCH models in order to test for normality of innovations,
        and it provides the pertaining large sample theory. Section \ref{sec_monte} presents a Monte Carlo study that compares the new tests with competing ones.
       The article concludes with discussions in Section \ref{sec_conclusion}.
       The majority of proofs has been postponed to Section \ref{sec_app}.

       Throughout  the paper, the letter $d$ stands for dimension, and both random and nonrandom vectors are understood as column vectors; the transpose of a vector $x$ will be denoted by $x^\top$; all random elements are supposed to be defined on the same probability space $(\Omega, \mathcal{A}, \mathbb{P})$; for any  matrix $A=(a_{kj})$, we will use the norm defined by $\| A \|=\sum_{k,j}|a_{k j}|$; if $A$ is a vector, $\| A \|$ denotes the Euclidean norm. If $A$ is a square matrix, $\|A\|_2$ stands for the spectral norm of $A$. The unit matrix of order $d$ will be denoted by I$_d$. Finally, ${\rm Re}(z)=a$ is the real part of a complex number $z=a +{\rm i}b$.

\section{Characterizations of multinormality}\label{sec_charac}
Let $X$ be a $d$-variate non-degenerate random vector with CF $\varphi_X(t) = \BE[\exp(\textrm{i}t^\top X)]$ and MGF $M_X(t) = \BE[\exp(t^\top X)]<\infty$,
 $t \in \RR^d$. We then have the following characterization.
\begin{prop} \label{char} The identity
\begin{equation} \label{CF1}
\varphi_X(t) \, M_X(t)-\e^{(\ii +1)t^\top\mu} \ =\ 0 \quad  \textrm{for each } t  \in \RR^d
\end{equation}
holds true for some $\mu \in \RR^d$ if,  and only if,  $X$ follows some normal distribution with mean $\mu$.
\end{prop}
\noindent {\sc Proof}.
The ``if'' part is trivial. To prove the converse implication, suppose that (\ref{CF1}) holds. Fix $a \in \RR^d$ and put $Y_a =a^\top (X-\mu)$. We have
$
\varphi_{Y_a}(t)=\e^{-\ii t a^\top \mu} \, \varphi_X(t a)$, $M_{Y_a}(t)=\e^{-t a^\top \mu} \, M_X(t a)$.
Hence, (\ref{CF1}) implies
$\varphi_{Y_a}(t) \, M_{Y_a}(t)-1 = 0$ for each $t\in \RR$.
In view of (\ref{volk}), it follows that $Y_a$ has some  zero-mean normal distribution. Moreover, the variance of $Y_a$ is equal to $a^\top \Sigma a$,
where $\Sigma$ denotes the covariance matrix of $X$. Hence, $a^\top X = Y_a + a^\top \mu$ has a normal distribution with mean  $a^\top \mu$
 and variance $a^\top \Sigma a$. Since $a$ was arbitrary, a well-known characterization of the multivariate normal distribution
  (see for instance \cite{rao:73}, \S8a.1--\S8a.2) yields the assertion.\bewend

The following result will be used in the construction of the test statistics and in order to  prove consistency of one of our tests.
\medskip
\begin{cor} \label{charzero}
Assume that $\BE X =0$ and put $R_X(t) := {\textrm{Re}} \, \varphi_X(t)$. Then
\begin{equation} \label{RCF}
R_X(t) \, M_X(t)-1 \ =\ 0 \quad  \textrm{for each } t  \in \RR^d
\end{equation}
holds true if,  and only if,  $X$ follows a zero-mean normal distribution.
\end{cor}
\noindent {\sc Proof}.
Since the ``if'' part is obviously true, suppose that (\ref{RCF}) holds. Replace $t$ with $-t$ to get
\begin{equation} \label{RCF1}
R_X(t) \, M_X(t)=1 \ =\ R_X(-t) \, M_X(-t) \quad  \textrm{for each } t  \in \RR^d.
\end{equation}
However, $R_X(t)=R_X(-t)$, and since, by (\ref{RCF}), $R_X(t)\neq 0$, (\ref{RCF1}) yields
\begin{equation} \label{RCF2}
 M_X(t)= M_X(-t) \quad  \textrm{for each } t  \in \RR^d.
\end{equation}
Noticing that $M_X(-t)=M_{-X}(t)$, (\ref{RCF2}) gives $M_X=M_{-X}$ which, by uniqueness of the MGF, shows that the law of $X$ is symmetric around zero,
 since $X$ and $-X$ have the same distribution. Hence $R_X = \varphi_X$, and Proposition \ref{char} completes the proof.  \bewend
%
%
%
%
\section{New tests for multivariate normality}\label{sec_teststat}
In this section, we assume that $X_1,X_2, \ldots, X_n, \ldots $ are i.i.d. copies of a $d$-dimensional random vector $X$,  the distribution of which
is assumed to be absolutely continuous with respect to $d$-dimensional Lebesgue measure. Writing N$_d(\mu,\Sigma)$ for the $d$-dimensional normal distribution
with mean vector $\mu$ and non-degenerate covariance matrix $\Sigma$ and ${\cal N}_d$ for the class of all non-degenerate $d$-variate normal distributions,
a classical problem is to test the null hypothesis
\[
H_0: \ {\rm{The \ law \ of}} \  X  \ {\rm{belongs\  to \ }} \ {\cal{N}}_d,
\]
against general alternatives. Since the class ${\cal N}_d$ is closed with respect to full rank affine transformations, any genuine test statistic $T_n = T_n(X_1,\ldots,X_n)$ based on $X_1,\ldots,X_n$ should also be invariant with respect to such transformations, i.e., we should have
\[
T_n(AX_1+b, \ldots, AX_n+b) \ = \ T_n(X_1,\ldots,X_n)
\]
for each nonsingular $d \times d$-matrix $A$ and each $b \in \RR^d$, see \cite{he:02} for an account on the importance of affine invariance in connection
with testing for multivariate normality. Writing $\overline{X}_n = n^{-1}\sum_{j=1}^n X_j$ for the sample mean and
 $S_n = n^{-1}\sum_{j=1}^n (X_j - \overline{X}_n)(X_j - \overline{X}_n)^\top$ for the sample covariance matrix of $X_1,\ldots,X_n$, a necessary and sufficient
 condition for a test statistic $T_n$ to be affine invariant is that it is based on the Mahanalobis angles and distances
 \begin{equation}\label{mahala}
 Y_{n,i}^\top Y_{n,j} \ = \ (X_i - \overline{X}_n)^\top S_n^{-1} (X_j - \overline{X}_n), \qquad i,j \in \{1,\ldots, n\},
 \end{equation}
 where
 \[
 Y_{n,j} = S_n^{-1/2}(X_j - \overline{X}_n), \qquad j=1,\ldots, n,
 \]
 are the so-called {\em scaled residuals} of $X_1,\ldots,X_n$, see \cite{he:02}.  Here, $S_n^{-1/2}$ denotes the unique symmetric square root of $S_n$ which,
 due to the absolute continuity of the distribution of $X$, exists with probability one if $n \ge d+1$, see \cite{ep:73}. The latter condition is tacitly assumed to hold in what follows.

Recall that in view of Corollary \ref{charzero}, we have
$
R_X(t) \, M_X(t) -1  =  0$ for each $t \in \RR^d$
if, and only if, the distribution of $X$ is centered normal. Since the scaled residuals $Y_{n,1}, \ldots, Y_{n,n}$ provide an empirical
standardization of $X_1,\ldots,X_n$, it seems tempting to introduce the {\em empirical cosine transform}
\begin{equation}\label{ECT}
R_n(t) \ := \ \frac{1}{n} \sum_{j=1}^n \cos\left(t^\top Y_{n,j}\right), \quad t \in \RR^d,
\end{equation}
 and the {\em empirical moment generating function}
\begin{equation}\label{EMGF}
M_n(t) \ := \ \frac{1}{n} \sum_{j=1}^n \exp\left(t^\top Y_{n,j}\right), \quad t \in \RR^d,
\end{equation}
of $Y_{n,1},\ldots,Y_{n,n}$ and to base a test of $H_0$ on the weighted $L^2$-statistic
\begin{equation}\label{teststat}
T_{n,\gamma} \ := \ \int_{\RR^d} U_n^2(t) \, w_\gamma(t) \, \textrm{d} t,
\end{equation}
where
\begin{equation}\label{unproc}
U_n(t) := \sqrt{n}\left(R_n(t)\, M_n(t) - 1 \right), \qquad t \in \RR^d,
\end{equation}
\begin{equation}\label{weightgamma}
w_\gamma(t) \ = \ \exp\left(-\gamma \|t\|^2\right),
\end{equation}
and $\gamma >0$ is some positive parameter.
In principle, we could replace $w_\gamma$ in (\ref{teststat}) with a more general weight function satisfying some general conditions.
The above special choice,
however,   leads to a test criterion with certain  extremely appealing features.
To this end, putting $Y_{jk}^\pm = Y_{n,j} \pm Y_{n,k}$, routine calculations give the representation
\begin{eqnarray}\label{Tn}
T_{n,\gamma} & = & \left(\frac{\pi}{\gamma}\right)^{d/2} \Bigg{\{}
\frac{1}{2n^3} \sum_{j,k,\ell,m =1}^n \Bigg[ \exp\left(\frac{\|Y^{+}_{jk}\|^2 - \|Y^{-}_{\ell m}\|^2}{4 \gamma} \right)
\cos\left(\frac{Y_{jk}^{{+\top}} Y^{-}_{\ell m}}{2\gamma } \right)\\ \nonumber
&+&\exp\left(\frac{\|Y^{+}_{jk}\|^2 - \|Y^{+}_{\ell m}\|^2}{4 \gamma} \right)
\cos\left(\frac{Y_{jk}^{{+\top}} Y^{+}_{\ell m}}{2\gamma } \right)\Bigg]\\ \nonumber
 \label{amen}
& & - \frac{2}{n} \sum_{j,k=1}^n \exp\left(\frac{\|Y_{n,j}\|^2- \|Y_{n,k}\|^2}{4\gamma} \right)\cos\left(\frac{Y_{n,j}^\top Y_{n,k}}{2\gamma } \right) + n \Bigg{\}},
\end{eqnarray}
which is amenable to computational purposes. Being a function of $ Y_{n,i}^\top Y_{n,j}$  figuring in (\ref{mahala}), the statistic
$T_{n,\gamma}$ is affine invariant. Rejection of $H_0$ is for large values of $T_{n,\gamma}$.



In view of the fact that $T_{n,\gamma}$ contains a fourfold sum which implies $O(n^4)$ additions to compute $T_{n,\gamma}$,  we also studied another statistic, which is
\begin{equation}\label{tntilde}
\widetilde{T}_{n,\gamma} \ :=  \int_{\RR^d} U_n(t) \, w_\gamma(t)  \, \textrm{d} t.
\end{equation}
A simple calculation shows that $\widetilde{T}_{n,\gamma}$ takes the form
\begin{equation} \label{tildetgamma}
\widetilde{T}_{n,\gamma} \ = \ \left( \frac{\pi}{\gamma}\right)^{d/2} \sqrt{n} \left( \frac{1}{n^2} \sum_{j,k=1}^n \exp \left(\frac{\|Y_{n,j}\|^2 - \|Y_{n,k}\|^2}{4 \gamma } \right)
 \cos \left(\frac{Y_{n,j}^\top Y_{n,k}}{2\gamma} \right) - 1\right)
\end{equation}
and hence is much faster to compute than $T_{n,\gamma}$.

We close this section by noting that tests for normality based on the empirical CF date back to \cite{ep:83}. For multivariate normality
the first contribution is by \cite{cs:86}, with later contributions by \cite{bh:88}, \cite{hz:90}, \cite{hw:97}, \cite{ep:99}, \cite{pu:05}, \cite{arc:07}, and \cite{te:11}.
 For review material on the empirical CF we refer to \cite{ush:99}. On the other hand, the approach based on the empirical MGF is certainly less popular,
 and it appears to include only a few entries. In the multivariate case we refer to the test of \cite{cw:89} for bivariate exponentiality testing,
 the extension of \cite{fll:98} of the $T_3$-plot of \cite{go:96} for testing normality  and to \cite{mh:10} for testing skew-normality.
 In the  univariate case we refer to \cite{me:07} and \cite{zg:10}, with the lacking theory of the latter test provided by \cite{hk:16}.
 Finally, for a recent general account of weighted $L^2$-statistics such as ours we refer to \cite{beh:16}.
%
%
%
%
%

\section{The case $\gamma \to \infty$}\label{secgaminfty}
In this section, we show that the statistic $T_{n,\gamma}$, after a suitable scaling, approaches a linear combination of two well-known
measures of multivariate skewness  as $\gamma \to \infty$. Likewise, $\widetilde{T}_{n,\gamma}$ is connected with a time-honored measure of multivariate kurtosis.

\begin{thm}\label{gammainf}
We have
\[
\lim_{\gamma \to \infty} \gamma^{3+d/2} \, \frac{96 T_{n,\gamma}}{n \pi^{d/2}} = 2b_{1,d} +  3\widetilde{b}_{1,d},
\]
where
\[
b_{1,d} = \frac{1}{n^2} \sum_{j,k=1}^n \left(Y_{n,j}^\top Y_{n,k}\right)^3, \qquad \widetilde{b}_{1,d} =  \frac{1}{n^2} \sum_{j,k=1}^n Y_{n,j}^\top Y_{n,k} \, \|Y_{n,j}\|^2 \, \|Y_{n,k}\|^2,
\]
are multivariate sample skewness in the sense of \cite{ma:70} and
\cite{mrs:93}, respectively.
\end{thm}
\noindent {\sc Proof}.  From (\ref{Tn}) and
\[
\exp(y) = 1 + y + \frac{y^2}{2} + \frac{y^3}{6} + O(y^4), \qquad \cos(y) = 1 - \frac{y^2}{2} + O(y^4)
\]
as $y \to 0$, the result follows by tedious but straightforward calculations, using the relations
$\sum_{j=1}^n Y_{n,j}  =  0$, $\sum_{j=1}^n \|Y_{n,j}\|^2 =  nd$, $\sum_{j,k=1}^n \|Y_{n,j}+ Y_{n,k}\|^2  =  2n^2d$,
\begin{eqnarray*}
\sum_{j,k=1}^n \|Y_{n,j}+ Y_{n,k}\|^4 & = & 2n^2 \left(\frac{1}{n} \sum_{j=1}^n \|Y_{n,j}\|^4 + d^2 + 2d\right), \\
\sum_{j,k=1}^n \|Y_{n,j}+ Y_{n,k}\|^4 Y_{n,j}^\top Y_{n,k}  & = &  8 \sum_{j,k=1}^n \left(Y_{n,j}^\top Y_{n,k}\right)^2 \|Y_{n,j}\|^2 + 4n^2 b_{1,d} + 2 n^2 \widetilde{b}_{1,d}
\end{eqnarray*}
as well as (writing tr$(D)$ for the trace of a square matrix $D$)
\begin{eqnarray*}
\sum_{j,k=1}^n \left(Y_{n,j}^\top Y_{n,k}\right)^2 \|Y_{n,j}\|^2  & = & \sum_{j,k=1}^n \textrm{tr}\left(Y_{n,j}^\top Y_{n,k}Y_{n,k}^\top Y_{n,j} \, \|Y_{n,j}\|^2\right) \\
& = & \textrm{tr}\left( \sum_{k=1}^n Y_{n,k}Y_{n,k}^\top \sum_{j=1}^n Y_{n,j}Y_{n,j}^\top \|Y_{n,j}\|^2\right)\\
& = & \textrm{tr}\left( n \textrm{I}_d \sum_{j=1}^n Y_{n,j}Y_{n,j}^\top \|Y_{n,j}\|^2\right)\\
& = & n \, \sum_{k=1}^n \textrm{tr}\left(Y_{n,j}^\top Y_{n,j} \|Y_{n,j}\|^2\right)\\
& = & n \sum_{j=1}^n \|Y_{n,j}\|^4.  \ \bewend
\end{eqnarray*}

\begin{rem} \label{remark.gammainf}
Theorem \ref{gammainf} parallels Theorem 2.1 of \cite{he:97}, who showed that the
BHEP statistic for testing for multivariate normality, after suitable rescaling, approaches $2b_{1,d} + 3 \widetilde{b}_{1,d}$ as a smoothing parameter
(called $\beta$ in that paper) tends to $0$. Since $\beta$ and $\gamma$ are related by $\beta = \gamma^{-1/2}$, this corresponds to letting $\gamma$ tend to infinity.
\end{rem}

\begin{thm}\label{gammainftilde}
We have
\[
\lim_{\gamma \to \infty} \gamma^{2+d/2} \, \frac{16  \widetilde{T}_{n,\gamma}}{\sqrt{n} \pi^{d/2}} = b_{2,d} -d(d+2),
\]
where
\[
b_{2,d} = \frac{1}{n} \sum_{j=1}^n \|Y_{n,j}\|^4
\]
is multivariate sample kurtosis in the sense of \cite{ma:70}.
\end{thm}
\noindent {\sc Proof}. From (\ref{tildetgamma}) we have
\begin{eqnarray*}
\gamma^{d/2} \frac{\widetilde{T}_{n,\gamma}}{\sqrt{n} \pi^{d/2}} & = &  \frac{1}{n^2}\sum_{i,j=1}^n \bigg{\{} \left( 1 + \frac{\|Y_{n,i}\|^2\! -\!  \|Y_{n,j}\|^2}{4 \gamma}
+ \frac{\left(\|Y_{n,i}\|^2\! - \! \|Y_{n,j}\|^2\right)^2}{32\gamma^2}\right) \\
& & \quad \times \left( 1- \frac{\left(Y_{n,i}^\top Y_{n,j}\right)^2}{8 \gamma^2}\right) + O\left(\gamma^{-3}\right) \bigg{\}} - 1
\end{eqnarray*}
as $\gamma \to \infty$.
Since $\sum_{i,j} (\|Y_{n,i}\|^2- \|Y_{n,j}\|^2) = 0$ and $\sum_{i,j} \left(Y_{n,i}^\top Y_{n,j}\right)^2 = n^2d$, the result follows by straightforward calculations. \ \bewend

\vspace*{2mm}
\noindent Notice that  $d(d+2)$ is the  value of theoretical multivariate kurtosis in the sense of \cite{ma:70} attained by any non-degenerate $d$-variate normal distribution. In this sense the limit statistic of Theorem \ref{gammainftilde} maybe termed ``excess kurtosis''.

%
%
%
%
%

\section{Asymptotic null distribution in the i.i.d. case}\label{sec_iid}
In this section we consider the case that $X_1,X_2, \ldots $ are i.i.d. $d$-dimensional random
vectors with some non-degenerate normal distribution.
Notice that the process $U_n$ given in (\ref{unproc}) is a random element of the separable Hilbert space
$\textrm{L}^2 := \textrm{L}^2(\RR^d,{\cal B}^d,w_\gamma(t)\textrm{d}t)$ of (equivalence classes of) measurable functions $f:\RR^d \rightarrow \RR$
that are square integrable with respect to the finite measure on the $\sigma$-field ${\cal B}^d$ of Borel sets of $\RR^d$ given by the weight function $w_\gamma$
defined in (\ref{weightgamma}).
The inner product and the resulting norm in $\textrm{L}^2$ will be denoted by
$
\langle f,g \rangle = \int_{\RR^d} f(t) \, g(t) \, w_\gamma(t) \, \textrm{d}t$, $\|f\|_{\textrm{L}^2} = \sqrt{\langle f,f \rangle}$,
respectively.

With this notation, the  test statistics $T_{n,\gamma}$  and $\widetilde{T}_{n,\gamma}$ given in (\ref{teststat}) and (\ref{tntilde}), respectively, take the form
\begin{equation}\label{unl2}
T_{n,\gamma} \ = \ \|U_n\|^2_{\textrm{L}^2}
\end{equation}
and
\begin{equation}\label{unl3}
\widetilde{T}_{n,\gamma} \ = \langle U_n,1 \rangle.
\end{equation}
Writing "$\vertk$" for convergence in distribution of random vectors and stochastic processes, the main result of this section
is as follows.

\smallskip
\begin{thm}{\rm{(}}Convergence of $U_n$ under $H_0${\rm{)}}\label{mainthm}\\
Suppose that $X$ has some non-degenerate $d$-variate normal distribution, and that $\gamma >1$ in (\ref{weightgamma}). Then there is a centered Gaussian random element $W$ of
$\LL$ having covariance kernel
\[
C(s,t)   =   {\rm{e}}^{s^\top t} + \frac{1}{2} \left({\rm{e}}^{s^\top t} + {\rm{e}}^{-s^\top t} \right) + 2 \cos\left(s^\top t\right)
- s^\top t  - 4
\]
so that $U_n \vertk W$ as $n \to \infty$.
\end{thm}

In view of \eqref{unl2} and  \eqref{unl3}, Theorem \ref{mainthm} and the Continuous Mapping Theorem immediately yield the following corollary.

\begin{cor} \label{cor1}
Under the assumptions of Theorem \ref{mainthm}, we have
\begin{enumerate}
\item[a)] $ T_{n, \gamma} \vertk \|W\|^2_{\LL} \ \ \textrm{as} \ n \to \infty$,
\item[b)] $\widetilde{T}_{n,\gamma} \vertk \textrm{N}(0, \sigma^2) \ \ \textrm{as} \ n \to \infty$,
with
$
\sigma^2 = \int_{\RR^d} \int_{\RR^d} C(s,t) {\rm{e}}^{-\gamma (\|s\|^2 + \|t\|^2)} \, \textrm{d} s \textrm{d}t.
$
\end{enumerate}
\end{cor}

\begin{rem} Routine calculations show that  $\sigma^2$ in Corollary \ref{cor1} (b) has the following expression:
$\sigma^2=2\pi^d(\gamma^2-0.25)^{-d/2}+2\pi^d(\gamma^2+0.25)^{-d/2}-4\pi^d\gamma^{-d}$.
\end{rem}

\begin{rem}
It is well-known that the distribution of $\|W\|^2_{\LL}$ is that of $\sum_{j=1}^\infty \lambda_jN_j^2$, where $\lambda_1, \lambda_2, \ldots $ are the
positive eigenvalues of the integral operator $f \mapsto Af$  on $\LL$ associated with the kernel $C$ given in Theorem \ref{mainthm},
i.e., $(Af)(t) \! = \! \int \! C(s,t) f(s) \exp(-\gamma \|t\|^2) \textrm{d}s$,
and $N_1,N_2, \ldots $ are i.i.d. standard normal random variables.
We did not succeed in obtaining explicit solutions of this equation.
However, since
\[
\mathbb E \|W\|^2_{\LL} = \int_{\RR^d} C(t,t) \, \exp\left( - \gamma \|t\|^2 \right) \, \textrm{d}t
\]
(see \cite{sw:86}, p. 213), tedious but straighforward manipulations of integrals yield
\begin{eqnarray} \nonumber
\mathbb E \|W\|^2_{\LL} & = & \frac{3}{2} \left(\frac{\pi}{\gamma-1}\right)^{d/2}+\frac{1}{2} \left(\frac{\pi}{\gamma+1}\right)^{d/2}-\left(4 +\frac{d}{2\gamma}\right)\left(\frac{\pi}{\gamma}\right)^{d/2}
\\
\label{mean}
& & + 2\left[\sqrt{\frac{\pi}{\gamma}} \frac{\cos\left(\frac{\tan^{-1}(1/\gamma)}{2}\right)}{\left(1+(1/\gamma^2)\right)^{1/4}}\right]^d \sum_{q=0}^{\lfloor \frac{d}{2}\rfloor } (-1)^q \: {d \choose 2q}
\left(\tan \left(\frac{\tan^{-1}(1/\gamma)}{2}\right)\right)^{2q},
\end{eqnarray}
\end{rem}
where $\lfloor \cdot \rfloor $ denotes the integer part of a given number.
%
%
%
%
%

\section{Consistency}\label{sec_consist}

The next result shows that the test for multivariate normality based on $T_{n,\gamma}$ is consistent against general alternatives.

\begin{thm}\label{thmconsist}
Suppose $X$ has some absolutely continuous distribution, and that $M_X(t) = \BE [\exp(t^\top X)] < \infty$, $t \in \RR^d$.
Then we have
\[ 
\liminf_{n \to \infty} \frac{T_{n,\gamma}}{n} \ \ge \  \int_{\RR^d} \left( R_X(t)M_X(t) -1 \right)^2 \, w_\gamma(t) \, \rm{d} t
\] 
almost surely.
\end{thm}

\noindent {\sc Proof}.  Because of affine invariance we may w.l.o.g. assume $\BE X =0$ and $\BE X X^\top = \textrm{I}_d$.
With $\Delta_{n,j}$ given in (\ref{delta1}), notice that
\[
\xi_n := \max_{j=1,\ldots,n} \|\Delta_{n,j}\| \le \|S_n^{-1/2} - \textrm{I}_d \|_2 \cdot \max_{j=1,\ldots,n} \|X_j\| + \|S_n^{-1/2}\|_2 \, \|\overline{X}_n\|
\]
which, by Theorem 5.2. of \cite{bn:63}, implies
\begin{equation}\label{maxabsch}
\lim_{n\to \infty} \xi_n = 0 \quad \PP\textrm{-a.s.}
\end{equation}
Fix $K>0$, and recall $R_n(t)$, $M_n(t)$, from (\ref{ECT}), (\ref{EMGF}), respectively. With  $R_n^\circ(t)$, $M_n^\circ(t)$ given in  (\ref{procrn0}),
(\ref{taylexp}) yields
\[
\max_{\|t\|\le K} \big{|} M_n(t) - M_n^\circ(t) \big{|} \le \max_{\|t \| \le K} M_n(t) \, \left(K \xi_n + \frac{1}{2} K^2 \xi_n^2  \textrm{e}^{K \xi_n}\right).
\]
By (\ref{maxabsch}) and the strong law of large numbers in the Banach space of continuous functions on $\{t \in \RR^d: \|t\| \le K\}$,
we have
$
\lim_{n \to \infty} \max_{\|t\|\le K} \big{|} M_n(t) - M_n^\circ(t) \big{|}   = 0
$
$\PP$-a.s. Since $\max_{\|t\| \le K}|R_n(t) - R_n^\circ(t)| \le K \xi_n$, (\ref{maxabsch}) implies
$
\lim_{n \to \infty} \max_{\|t\|\le K} \big{|} R_n(t) - R_n^\circ(t) \big{|}   = 0$  $\PP\textrm{-a.s.}$
In view of
\begin{eqnarray*}
R_n(t)M_n(t) - 1 & = & (R_n(t) - R_n^\circ(t))(M_n(t) - M_n^\circ(t)) + M_n^\circ(t)(R_n(t) - R_n^\circ(t))\\
& & + R_n^\circ (t)(M_n(t) - M_n^\circ(t)) + R_n^\circ(t)M_n^\circ(t)-1,
\end{eqnarray*}
the strong law of large numbers in Banach spaces and Fatou's lemma yield
\[ 
\liminf_{n \to \infty} \frac{T_{n,\gamma}}{n} \ \ge \  \int_{\{t:\|t\|\le K\}} \left( R_X(t)M_X(t) -1 \right)^2 \, w_\gamma(t) \, \mbox{d} t
\] 
almost surely. Since $K$ is arbitrary, the assertion follows. \ \bewend


In view of Corollary \ref{charzero} and Theorem \ref{thmconsist}, the test that rejects the null hypothesis ${ {H}}_0$ for large values of $T_{n,\gamma}$
is consistent against each alternative distribution the MGF of which exists on $\RR^d$. We conjecture that this test is consistent against {\em any} alternative distribution.
However, in view of the reasoning of \cite{cs:89}, the behavior of $T_{n,\gamma}$ against heavy-tailed alternatives is a non-trivial problem.


%
%
%
%
%

\section{Testing for normality in GARCH models}\label{sec_garch}
Consider the multivariate GARCH (MGARCH) model
\begin{equation} \label{GARCH}
X_j=\Sigma_j^{1/2}(\theta)\varepsilon_j, \ j\in \mathbb{Z},
\end{equation}
where $\theta \in \Theta \subseteq \mathbb R^v$, is an $v$-dimensional vector of  unknown parameters.
The unobservable random errors $\varepsilon_j$, $j \in \mathbb{Z}$, also referred to as {\em innovations},  are i.i.d. copies
 of a $d$-dimensional random vector $\varepsilon$, which is assumed to have mean zero and unit covariance matrix. Moreover,
$\Sigma_j(\theta)=\Sigma(X_{j-1},X_{j-2},\dots;\theta)$ is a $d\times d$ symmetric and positive definite matrix,
which is the conditional variance matrix of $X_j$, given the past information.
On the basis of the observations $X_j, \ j=1,\ldots,n$, driven by equation (\ref{GARCH}), we wish to test the null
hypothesis
\[
{H}_{0,G}: \ {\rm{The \ law \ of}} \  \varepsilon  \ {\rm{is}} \ \textrm{N}_d(0,{ {\rm I}}_d).
\]
Notice that ${{H}}_{0,G}$ is equivalent to the hypothesis that, conditionally on $\{X_{j-1},X_{j-2},\ldots\}$, the law of $X_j$  is N$_d(0,\Sigma_j(\theta))$,
for some $\theta \in \Theta$. The difference with the i.i.d. setting is that, on the one hand,
the innovations in (\ref{GARCH}) are already assumed to be centered at zero. On the other hand however,
the covariance matrix $\Sigma_j(\theta)$ of $X_j$, unlike the i.i.d. case, here is allowed to be time-varying in a way that depends on the unknown parameter $\theta$ as well as on past observations.

For univariate GARCH models \cite{klm:12}, \cite{llp:14}, \cite{gr:14}  and \cite{dmb:17}
suggested specification tests for the innovation distribution. However, with the exception of \cite{bc:08}, corresponding tests
 are still scarce in the multivariate case. We now take some time to emphasize the importance of testing the null hypothesis of normality in GARCH models.
 First, notice that acceptance of the null hypothesis  ${{H}}_{0,G}$  implies the validity of the classical Gaussian MGARCH model,
 which has been a benchmark for modelling certain economic and financial quantities. Although even today the normal distribution is
 most commonly used in applications, since the time of \cite{ma:63} and \cite{fa:65} there is empirical evidence that,
for example, the distribution of financial variables is heavy-tailed, even after filtering the volatility clustering phenomenon produced by model  (\ref{GARCH}).
In order to capture this empirical so-called stylized fact, several authors suggested alternative innovation distributions, such as
 the (multivariate) Student-t distribution (\cite{tsa:06}, \cite{bc:08}, and \cite{dv:10}), the stable distribution (\cite{bo:12}, \cite{og:13}, and \cite{fm:16}),
  and the Laplace distribution (\cite{tz:07}, \cite{cgl:12}). In this connection, it is well-known that having erroneously
  accepted the normality assumption for the GARCH-residuals resulting from the estimation
of model (\ref{GARCH}), leads to incorrect inferential procedures, such as assessment of standard risk measures like the value at risk (VaR); see for instance \cite{sp:16}.
The preceding discussion provides the ground on the basis of which the null hypothesis ${{H}}_{0,G}$ could be considered as highly relevant, particularly in statistical modelling with a view towards financial applications.

Notice that although  ${{H}}_{0,G}$  is about the distribution of $\varepsilon$, the innovations themselves are unobservable in the context of model (\ref{GARCH}). Hence any decision regarding the  null hypothesis  ${{H}}_{0,G}$ should be based on residuals
\begin{equation} 
{\widetilde {\varepsilon}}_j(\widehat{\theta}_n)={\widetilde{ \Sigma}}_j^{-1/2}(\widehat{\theta}_n)X_j, \ j=1,\ldots,n.
\end{equation}
Here, the matrix $\widetilde{ \Sigma}_j(\cdot)$, apart from a suitable estimator $\widehat{\theta}_n$ that will be detailed later, also depends on
specific initial values $\{\widetilde X_j, \ j\leq 0\}$ of $X_j$ which, under certain conditions, are asymptotically irrelevant. Let $U_n^G$ be defined  as $U_n$ in \eqref{unproc} by replacing $Y_{n,j}$ with ${\widetilde {\varepsilon}}_j(\widehat{\theta}_n), \ j=1,\ldots,n$.
The value of the test statistics $T_{n, \gamma}^G$  and $\widetilde{T}_{n,\gamma}^G$ are  defined as in  (\ref{teststat}) and (\ref{tntilde}), respectively,  with $U_n$ changed for $U_n^G$.

In order to derive the asymptotic null distribution of $U_n^G$ we will make the following assumptions (A.1)--(A.6), which will be commented on at the end of this section for specific instances of MGARCH models.
  In the sequel, $C>0$ and $\varrho$, $0< \varrho < 1$,  denote generic constants,
 the values of which may vary across the text, and $\theta_0$ stands for the true value of $\theta$.

 \begin{enumerate}
\item[(A.1)] The estimator $\widehat{\theta}_n$ satisfies
$
\sqrt{n}(\widehat{\theta}_n-\theta_0)={O}_{\mathbb{P}}(1).
$
%

\item[(A.2)] $\sup_{\theta\in\Theta}\left\|\widetilde{\Sigma}^{-1/2}_{j}(\theta)\right\|\leq C, \quad \sup_{\theta\in\Theta}\left\|\Sigma^{-1/2}_{j}(\theta)\right\|\leq C\quad \mbox{a.s.}$

\item[(A.3)] $\sup_{\theta\in\Theta}\|\Sigma^{1/2}_{j}(\theta)-\widetilde{\Sigma}^{1/2}_{j}(\theta)\| \leq C\varrho^j.$

\item[(A.4)] $\mathbb E \left\|X_j\right\|^\varsigma<\infty$ and $ \mathbb E\left\|\Sigma^{1/2}_{j}(\theta_0)\right\|^\varsigma<\infty$  for some $\varsigma>0$.

\item[(A.5)] For each sequence $x_1,x_2,\dots$ of vectors of $\mathbb{R}^d$, the function $\theta\mapsto	
	\Sigma^{1/2}(x_1,x_2,\dots;\theta)$ admits continuous second-order derivatives.

\item[(A.6)] For some neighborhood $V(\theta_0)$ of  $\theta_0$, there exist $p> 1$, $q> 2$ and $r> 1$ so that $2p^{-1}+2r^{-1}=1$ and $4q^{-1}+2r^{-1}=1$, and
    \begin{eqnarray*}
&& \mathbb E \sup_{\theta\in V(\Theta)}\left\|\sum_{k,\ell=1}^v\Sigma_j^{-1/2}(\theta)\frac{\partial^2\Sigma^{1/2}_j(\theta)}{\partial \theta_k\partial \theta_\ell}\right\|^p<\infty,\\
&& \mathbb E \sup_{\theta\in V(\Theta)}\left\|\sum_{k=1}^v\Sigma_j^{-1/2}(\theta)\frac{\partial\Sigma^{1/2}_j(\theta)}{\partial \theta_k}\right\|^q<\infty,\\
&& \mathbb E \sup_{\theta\in V(\Theta)}\left\|\Sigma_j^{1/2 }(\theta_0)\Sigma_j^{-1/2}(\theta)\right\|^r<\infty.
	\end{eqnarray*}
\end{enumerate}

The next result gives the asymptotic null distribution of $U_n^G$.

\begin{thm}{\rm{(}}Convergence of $U_n^G$ under $H_{0,G}${\rm{)}}\label{mainthmG}\\
Let $\{X_j\}$ be a strictly stationary process satisfying \eqref{GARCH}, with $X_j$ being measurable with respect to the sigma-field generated by $\{\varepsilon_u,u\leq j\}$. Assume  that
    (A.1)--(A.6) hold. Then under the null hypothesis ${{H}}_{0,G}$,  there is a centered Gaussian random element $W_G$ of
$\LL$ having covariance kernel
\[
C_G(s,t)   =   {\rm{e}}^{s^\top t} + \frac{1}{2} \left({\rm{e}}^{s^\top t} + {\rm{e}}^{-s^\top t} \right) + 2 \cos\left(s^\top t\right) - 4
\]
so that $U_n^G \vertk W_G$ as $n \to \infty$.
\end{thm}

From Theorem \ref{mainthmG} and the Continuous Mapping Theorem we have the following corollary.

\begin{cor} \label{cor2}
Under the assumptions of Theorem \ref{mainthmG}, and for $\gamma >1$ in (\ref{weightgamma}) we have
\begin{enumerate}
\item[a)] $ T_{n, \gamma}^G \vertk \|W_G\|^2_{\LL} \ \ \textrm{as} \ n \to \infty$,
\item[b)] $\widetilde{T}_{n,\gamma}^G \vertk \textrm{N}(0, \sigma^2_G) \ \ \textrm{as} \ n \to \infty$,
where
$
\sigma^2_G = \int_{\RR^d} \int_{\RR^d} C_G(s,t) {\rm{e}}^{-\gamma (\|s\|^2 + \|t\|^2)} \, \textrm{d} s \textrm{d}t.
$
\end{enumerate}
\end{cor}

\begin{rem} \label{remG}
Notice that the covariance kernels in Theorems \ref{mainthm} and \ref{mainthmG} only differ by the term $s^\top t$. This difference is due to the fact
that for the MGARCH case the innovations (and thus the data) are assumed to have zero mean. Hence, there is no centering involved.
 Second, since $\int_{\RR^d} \int_{\RR^d} s^\top t {\rm{e}}^{-\gamma (\|s\|^2 + \|t\|^2)} \, \textrm{d} s \textrm{d}t=0$, it follows that $\sigma^2=\sigma^2_G$,
and thus $\widetilde{T}_{n,\gamma}$ and $\widetilde{T}_{n,\gamma}^G$ both have the same asymptotic null distribution. Perhaps the most surprising fact is that
  the limit null distribution of $U_n^G$ does not depend on the estimation of the parameter $\theta$, in contrast to the GARCH version of the normality test of
  \cite{hw:97}, studied in \cite{klm:12}, \cite{jg:14} and \cite{jp:17} for the usual (linear) univariate GARCH model.
   We  underline that these observations refer to the asymptotic null distribution alone and that for finite sample sizes,
    the null distribution of both $ T_{n, \gamma}^G $ and $\widetilde{T}_{n,\gamma}^G$ do depend  on the estimation of the parameter $\theta$ as well as on the true value of this parameter.
\end{rem}

\begin{rem}
\cite{cl:03},  \cite{lmc:03}, \cite{bw:09} and  \cite{fz:12}, among others,
have shown that mild regularity conditions guarantee that the quasi maximum likelihood estimator (QMLE), defined by
\[ 
\widehat{\theta}_n=\argmax_{\theta \in \Theta} {\cal{L}}_n(\theta),
\] 
\[
{\cal{L}}_n(\theta)=-\frac{1}{2}\sum_{j=1}^n {\widetilde {\ell}_j}(\theta),\quad{\widetilde {\ell}_j}(\theta )=X_j^\top
\widetilde{ \Sigma}_j^{-1}(\theta)X_j
+\log\left|\widetilde{ \Sigma}_j (\theta)\right|,
\]
 satisfies (A.1) for general MGARCH,
or for particular specifications.

There are many MGARCH parametrizations for the matrix $\Sigma_j(\theta)$, see, e.g., \cite{fz:10}.
One of the most widely used MGARCH models is the Constant Conditional Correlation (CCC) model proposed by \cite{bo:90} and  extended by \cite{jt:98}.
That model decomposes the conditional covariance matrix (\ref{GARCH}) into conditional standard
deviations and a conditional correlation matrix, according to
$
\Sigma_j(\theta_0)={ {D}}_j(\theta_0) R_0 {{D}}_j(\theta_0).
$
Here, ${{D}}_j(\theta_0)$ and $R_0$ are $d\times d$ matrices, with $R_0$ being a correlation matrix and  ${{D}}_j(\theta_0)$ is a diagonal
 matrix so that ${\sigma}_j(\theta)=\mbox{diag}\left\{D^2_j(\theta)\right\}$ with
\[
{{\sigma}}_j(\theta)={b}+\sum_{k=1}^p {{B}}_k X_{j-k}^{(2)}+\sum_{k=1}^q {{\Gamma}}_k {{\sigma}}_{j-k}(\theta)
\]
and $X_{j}^{(2)}= X_{j} \odot X_{j}$, where $\odot$ denotes the Hadamard product, that is, the element by element product.
Moreover,
${{b}}$ is a vector of dimension $d$ and has positive elements, while  $\{B_k\}_{k=1}^p$ and $\{{{\Gamma}}_k\}_{k=1}^q$ are  $d\times d$ matrices  with non-negative elements.
Under certain weak assumptions, the QMLE for the parameters in this model satisfies (A.1), and (A.2)--(A.6) also hold (see \cite{fz:10}  and \cite{fjm:17}).
\end{rem}
\begin{rem}
Similar results for the consistency against alternatives to those stated in the i.i.d case can be given now. To save space we omit them.
\end{rem}
%
%
%
%
%
\section{Monte Carlo results}\label{sec_monte}
This section describes and summarizes the results of some simulation experiments. All computations have
been performed using programs written in the R language.

\subsection{Numerical experiments for i.i.d. data}
Upper quantiles of the null distribution of $T_{n,\gamma}$ have been approximated by  gene\-ra\-ting 10,000 samples  from a law N$_d(0,\textrm{I}_d)$. Table  \ref{tabla2}  displays some critical values. Looking at this table we see that the speed of convergence to the asymptotic null distribution depends on the data dimension and the value of $\gamma$.
\begin{table} 
\centering
\caption{Critical points for $(\gamma/\pi)^{d/2}T_{n,\gamma}$.} \label{tabla2}
\begin{tabular}{|ccc|rrrrrrr|}
   \hline
    &     &          & \multicolumn{7}{c|}{$\gamma$}\\ \cline{4-10}
$d$ & $n$ & $\alpha$ & 1.2 & 1.3 & 1.4 & 1.5 &  2.0 &  2.5 &  3.0  \\ \hline
2  & 20 & 0.05 &      7.95 &    5.62 &   3.94 &   2.71 &   0.80 &  0.31 &  0.15 \vspace{-3pt}\\
   &    & 0.10 &      4.38 &    2.94 &   2.08 &   1.53 &   0.46 &  0.19 &  0.10 \vspace{-3pt}\\
   & 50 & 0.05 &     12.62 &    7.19 &   4.82 &   3.74 &   1.06 &  0.43 &  0.21 \vspace{-3pt}\\
   &    & 0.10 &      6.39 &    4.06 &   2.81 &   2.00 &   0.62 &  0.26 &  0.13 \vspace{-3pt}\\  
  & 100 & 0.05 &     14.37 &    8.11 &   5.22 &   3.60 &   1.04 &  0.44 &  0.22 \vspace{-3pt}\\
   &    & 0.10 &      6.97 &    4.31 &   2.98 &   2.15 &   0.67 &  0.29 &  0.14  \vspace{-3pt}\\  \hline
3  & 20 & 0.05 &     22.71 &   13.73 &   9.20 &   6.75 &   1.84 &  0.70 &  0.33 \vspace{-3pt}\\
   &    & 0.10 &     13.37 &    8.44 &   5.80 &   4.29 &   1.22 &  0.48 &  0.23 \vspace{-3pt}\\
   & 50 & 0.05 &     50.43 &   24.12 &  14.41 &   9.17 &   2.36 &  0.93 &  0.45 \vspace{-3pt}\\
   &    & 0.10 &     23.51 &   12.45 &   7.92 &   5.63 &   1.54 &  0.63 &  0.31 \vspace{-3pt}\\
  & 100 & 0.05 &     65.97 &   28.84 &  16.20 &   9.81 &   2.35 &  0.95 &  0.46 \vspace{-3pt}\\
   &    & 0.10 &     28.58 &   15.06 &   9.07 &   6.02 &   1.61 &  0.66 &  0.33 \vspace{-3pt}\\ \hline
5  & 20 & 0.05 &    124.67 &   63.48 &  37.98 &  25.23 &   5.99 &  2.18 &  0.98 \vspace{-3pt}\\
   &    & 0.10 &     74.76 &   41.77 &  26.13 &  17.49 &   4.44 &  1.67 &  0.77 \vspace{-3pt}\\
   & 50 & 0.05 &    670.14 &  222.95 & 100.01 &  54.12 &   8.46 &  3.05 &  1.41 \vspace{-3pt}\\
   &    & 0.10 &    264.18 &  105.49 &  50.71 &  29.23 &   5.72 &  2.17 &  1.03 \vspace{-3pt}\\
  & 100 & 0.05 &   1501.74 &  428.77 & 157.74 &  73.21 &   9.13 &  3.22 &  1.48 \vspace{-3pt}\\
   &    & 0.10 &    514.14 &  170.50 &  71.73 &  36.88 &   5.98 &  2.27 &  1.09 \vspace{-3pt}\\ \hline
\end{tabular}
\end{table}

A natural competitor of the test based on $T_{n,\gamma}$ is the CF-based test studied in  \cite{hw:97} (HW-test), since this test is simple to compute
as well as  affine invariant and has been shown to have good power performance vis--\'a-vis competitors.
Since, according to \cite{ja:00}, the global power function of any nonparametric test is flat on balls of alternatives
  except for alternatives coming from a finite-dimensional subspace, each test has a high power only against  a specific set of alternatives.
So, we carried out an extensive simulation study to compare their powers against a wide range of alternatives, with the aim of detecting those for which the test
based on $T_{n,\gamma}$ is more powerful than the HW-test. As expected from Theorem \ref{gammainf} and Remark \ref{remark.gammainf},  for large $\gamma$ and small
$\beta$ both tests behave very closely. For non-heavy-tailed distributions, we observed that the power of the proposed test is either similar
or a bit less than that of the HW-test; for heavy-tailed distributions, however, the new test outperforms the HW-test. This observation can be appreciated by looking at
Table \ref{tabla3}, which displays the empirical power calculated by generating 1,000 samples (in each case), for the significance level $\alpha=0.05$,
from the following heavy-tailed alternatives: ($LA$) the multivariate symmetric Laplace distribution (\cite{kp:01});
($GN_{\theta}$) the multivariate $\beta$-generalized distribution (\cite{gk:73}), that coincides with the normal distribution for $\theta=2$ and has heavy tails for $0<\theta<2$;
($ASE_{\theta}$) the $\theta$-stable  and  elliptically-contoured distribution;
($T_{\theta}$) multivariate student-$t$ with $\theta$ degrees of freedom.

\begin{table} 
\vspace*{-2cm}
\centering
\caption{Percentage of rejection for nominal level $\alpha=0.05$ and $n=50$.} \label{tabla3}
\footnotesize
\begin{tabular}{|cc|rrrrr|rrr|}
\hline
     & &   \multicolumn{5}{c|}{Test based on $T_{n,\gamma}$} & \multicolumn{3}{c|}{HW-test }\\ \cline{3-10}
     & &   \multicolumn{5}{c|}{$\gamma$} & \multicolumn{3}{c|}{$\beta$}\\  \hline

   & $d$ &   1.3&  1.4 &  1.5 & 2.0 & 2.5 &  0.1 & 0.5 & 1.0\\
   \hline
$LA$
& 2 & 83.6 & 83.8 & 82.8 & 82.0 & 81.1 & 59.5 & 75.3 & 82.5 \vspace{-3pt}\\
& 3 & 91.9 & 92.8 & 93.8 & 93.7 & 92.8 & 77.0 & 89.5 & 94.2 \vspace{-3pt}\\
& 5 & 96.9 & 97.6 & 97.9 & 99.0 & 98.8 & 96.0 & 98.6 & 99.4 \vspace{-3pt}\\
\hline$T_5$
& 2 & 66.7 & 67.2 & 65.9 & 65.5 & 65.0 & 52.9 & 56.2 & 50.0 \vspace{-3pt}\\
& 3 & 75.5 & 76.6 & 78.3 & 78.7 & 78.0 & 67.5 & 69.9 & 59.6 \vspace{-3pt}\\
& 5 & 87.1 & 88.2 & 88.6 & 91.0 & 91.5 & 86.0 & 88.0 & 78.1 \vspace{-3pt}\\ \hline
$T_{10}$
& 2 & 31.4 & 33.7 & 34.1 & 32.4 & 33.2 & 25.8 & 23.9 & 18.5 \vspace{-3pt}\\
& 3 & 38.4 & 40.0 & 40.0 & 41.3 & 41.3 & 33.7 & 31.7 & 22.3 \vspace{-3pt}\\
& 5 & 47.9 & 50.1 & 50.6 & 51.6 & 55.4 & 51.1 & 45.3 & 27.4 \vspace{-3pt}\\
\hline
$GN_{1}$
& 2 & 87.7 & 87.8 & 85.7 & 83.5 & 82.1 & 59.6 & 74.5 & 87.5 \vspace{-3pt}\\
& 3 & 91.3 & 92.6 & 93.1 & 92.1 & 91.4 & 72.0 & 88.6 & 95.9 \vspace{-3pt}\\
& 5 & 92.2 & 93.4 & 94.7 & 96.5 & 96.3 & 86.4 & 95.2 & 97.9 \vspace{-3pt}\\
\hline
$GN_{1.5}$
& 2 & 51.8 & 52.2 & 49.4 & 46.9 & 45.3 & 29.0 & 35.9 & 49.3 \vspace{-3pt}\\
& 3 & 56.2 & 58.5 & 60.7 & 59.0 & 56.3 & 38.9 & 48.5 & 66.6 \vspace{-3pt}\\
& 5 & 50.1 & 53.4 & 55.1 & 61.5 & 59.9 & 46.1 & 55.7 & 66.5 \vspace{-3pt}\\
\hline
$ASE_{1.75}$
& 2 & 75.6 & 75.8 & 75.5 & 75.4 & 75.8 & 68.5 & 69.2 & 60.9 \vspace{-3pt}\\
& 3 & 82.6 & 82.6 & 83.0 & 83.2 & 83.6 & 79.2 & 78.5 & 68.2 \vspace{-3pt}\\
& 5 & 89.0 & 89.1 & 89.4 & 90.6 & 90.7 & 88.1 & 86.0 & 74.0 \vspace{-3pt}\\
\hline
$ASE_{1.85}$
& 2 & 56.1 & 56.1 & 55.5 & 56.2 & 56.3 & 50.3 & 48.8 & 39.5 \vspace{-3pt}\\
& 3 & 63.2 & 63.7 & 64.4 & 64.1 & 63.6 & 57.5 & 55.5 & 42.8 \vspace{-3pt}\\
& 5 & 76.0 & 76.5 & 76.6 & 78.2 & 78.1 & 72.6 & 68.1 & 48.9 \vspace{-3pt}\\
\hline
\end{tabular}
\end{table}

As shown in Corollary \ref{cor1},  the test statistic  $\widetilde{T}_{n,\gamma}$ is asymptotically normal with zero mean under the null hypothesis. We carried out some simulations to assess the  normal approximation to the null distribution of $\widetilde{T}_{n,\gamma}$,  and observed that it requires very large values of $n$, which depend on the value of $\gamma$ and the data dimension.
 Table \ref{tabla5} exhibits the empirical power, for significance level $\alpha=0.05$, against some heavy-tailed  and  light-tailed distributions. We write $PII_{a}$ for the Pearson type II distribution with parameter $a$, and $U(0,1)^d$ for the uniform distribution on the $d$--dimensional cube. In Table \ref{tabla5} the one-sided test  that rejects $H_0$ for large values of   $\widetilde{T}_{n,\gamma}$ is codified as ``one'' while the two-sided test that rejects $H_0$ for large values of   $|\widetilde{T}_{n,\gamma}|$ is codified as ``two'' with critical points calculated by simulation. Such results were calculated by generating 1,000 samples of size $n=100$ in each case.
For very heavy-tailed distributions, the new test is clearly more powerful than the HW-test. In these cases, the one-sided test  gives slightly better results than the two-sided one.
For distribution with very light tails, the one-sided test fails (same behaviour as  the ``quadratic" statistic), but the two-sided test is more powerful than the HW-test. These numerical results are in agreement with the statement of Theorem \ref{gammainftilde}, which asserts that $\widetilde{T}_{n,\gamma}$ is close to a sample kurtosis measure.

\begin{table} 
\centering
\caption{Percentage of rejection for nominal level $\alpha=0.05$ and $n=100$.} \label{tabla5}
\footnotesize
\begin{tabular}{|c|c|lrrrrr|rrr|}
\hline
 & & \multicolumn{6}{c|}{Test based on $\widetilde{T}_{n,\gamma}$} & \multicolumn{3}{c|}{HW-test}\\
 \cline{3-11}
 & & \multicolumn{6}{c|}{$\gamma$} & \multicolumn{3}{c|}{$\beta$}\\
 \cline{3-11}
     & $d$    &     &    1.3 & 1.4 & 1.5 & 2.0 & 2.5  & 0.1 & 0.5 & 1.0 \vspace{-3pt}\\ \hline
$T_{10}$ & 2  & one &    57.0 & 57.0 & 57.1 & 57.7 & 57.9 & 33.3 & 34.8 & 27.1\vspace{-3pt}\\
     &        & two &    46.1 & 46.1 & 46.2 & 46.4 & 46.4 &   &  &\vspace{-3pt}\\
     & 3      & one &    74.9 & 74.7 & 74.6 & 74.4 & 74.1 & 44.7 & 47.4 & 38.9 \vspace{-3pt}\\
     &        & two &    64.7 & 64.7 & 64.8 & 65.5 & 65.3 &   &  &\vspace{-3pt}\\
     & 5      & one &    92.9 & 92.7 & 92.8 & 92.7 & 92.8 & 68.4 & 69.4 & 54.4 \vspace{-3pt}\\
     &        & two &    88.0 & 88.3 & 88.3 & 87.9 & 87.8 &   & &  \vspace{-3pt}\\ \hline
$ASE_{1.95}$ & 2 & one & 38.1 & 38.6 & 38.8 & 39.1 & 39.1 & 32.1 & 29.0 & 21.0 \vspace{-3pt}\\
     &           & two & 33.9 & 34.1 & 34.2 & 34.4 & 34.5 & &   &\vspace{-3pt}\\
     & 3         & one & 45.8 & 46.1 & 46.3 & 46.6 & 46.5 & 39.9 & 35.9 & 24.3 \vspace{-3pt}\\
     &           & two & 41.5 & 41.6 & 41.6 & 41.9 & 41.8 & &   &\vspace{-3pt}\\
     & 5         & one & 56.1 & 56.0 & 55.8 & 56.2 & 56.2 & 50.2 & 42.6 & 27.9 \vspace{-3pt}\\
     &           & two & 52.4 & 52.4 & 52.3 & 52.1 & 518 & &  & \vspace{-3pt}\\ \hline
     $U(0,1)^d$   & 2 & one &  0.0 & 0.0 & 0.0 & 0.0 & 0.0 & 2.0 & 28.3 & 97.4 \vspace{-3pt}\\
     &           & two & 99.9 & 100.0 & 100.0 & 100.0 & 100.0 & &   &\vspace{-3pt}\\
     & 3         & one & 0.0 & 0.0 & 0.0 & 0.0 & 0.0 & 0.0 & 27.0 & 98.0 \vspace{-3pt}\\
     &           & two & 100.0 & 100.0 & 100.0 & 100.0 & 100.0 & &  &\vspace{-3pt}\\
     & 5         & one & 0.0 & 0.0 & 0.0 & 0.0 & 0.0 & 0.0 & 22.4 & 95.8 \vspace{-3pt}\\
     &           & two & 100.0 & 100.0 & 100.0 & 100.0 & 100.0 & & &  \vspace{-3pt}\\ \hline
$PII_{4}$    & 2 & one & 0.0 & 0.0 & 0.0 & 0.0 & 0.0 & 3.0 & 1.7 & 11.6 \vspace{-3pt}\\
     &           & two & 47.8 & 48.0 & 48.2 & 48.1 & 48.6 & &  & \vspace{-3pt}\\
     & 3         & one & 0.0 & 0.0 & 0.0 & 0.0 & 0.0 & 0.2 & 1.7 & 17.9 \vspace{-3pt}\\
     &           & two & 72.7 & 73.2 & 73.9 & 74.8 & 75.0 & &  & \vspace{-3pt}\\
     & 5         & one & 0.0 & 0.0 & 0.0 & 0.0 & 0.0 & 0.0 & 2.7 & 32.6 \vspace{-3pt}\\
     &           & two & 93.9 & 94.5 & 94.8 & 95.3 & 95.6 & & & \vspace{-3pt}\\ \hline

\end{tabular}
\end{table}

\subsection{Numerical experiments for GARCH models}
Since usual practical applications of MGARCH models involve rather large sample sizes, this subsection studies the finite sample size behavior of the test based on $\widetilde{T}_{n,\gamma}^G$. With this aim, we must first specify a form for $\Sigma_j(\theta)$. In our simulations we considered a bivariate CCC--GARCH(1,1) model with
\[ {b}=\left(\begin{array}{c}0.1\\ 0.1\end{array}\right), \
 {{B}}_1= \left(\begin{array}{cc}0.3 & 0.1\\ 0.1 & 0.2\end{array}\right), \
 {{\Gamma}}_1= \left(\begin{array}{cc}0.2 & 0.1\\ 0.01 & 0.3 \end{array}\right), \
 {{R}}=\left(\begin{array}{cc}1 & r\\ r & 1 \end{array}\right),
\]
for $r=0, 0.3$, and a trivariate bivariate CCC--GARCH(1,1) model with  ${b}=(0.1, 0.1, 0.1)'$,
\[
 {B}_1= \left(\begin{array}{ccc}0.3 & 0.1 & 0.1\\ 0.1 & 0.2 & 0.1 \\ 0.1 & 0.1 & 0.1\end{array}\right), \
 {\Gamma}_1= \left(\begin{array}{ccc}0.2 & 0.1 & 0.01 \\ 0.01 & 0.3 & 0.1\\  0.01 & 0.01 & 0.1 \end{array}\right), \
 {R}=\left(\begin{array}{ccc}1 & r & r\\ r & 1 & r\\ r & r & 1 \end{array}\right)
\]
and $r$ as before.

As shown in Corollary \ref{cor2},  the test statistic  $\widetilde{T}_{n,\gamma}^G$ is asymptotically normal with zero mean under the null hypothesis.
We have computed the actual level for the one-sided test that rejects $H_{0,G}$ if $\widetilde{T}_{n,\gamma}^G/\sigma>u_{1-\alpha}$ and for the two-sided test rejecting $H_{0,G}$ when $|\widetilde{T}_{n,\gamma}^G/\sigma|>u_{1-\alpha/2}$ by generating 5,000 samples from the above bivariate CCC--GARCH(1,1) model.
The results are not reported in order to save space. We nevertheless have the following observations:
The quality of the  normal approximation strongly depends on the model parameter values, being conservative for $r=0.0$ and a bit liberal for $r=0.3$.
As in the i.i.d.-case, rather large sample sizes are required for the normal approximation to work. In addition, the normal approximation gives actual levels closer
to the nominal ones for the two-sided test. Nevertheless, as discussed before, there is empirical evidence that in typical applications the innovations
have heavy-tailed distributions, and we have learnt from the i.i.d. setting that in such a case, the one-sided test is more powerful than the two-sided one.
So, for moderate sample sizes, we should resort to another null distribution approximation.

In view of the remarks at the end of Remark \ref{remG}, and in contrast to the i.i.d.-case, we cannot calculate critical points for each $n$ by simulation since, in practice, the values of the true parameters involved in the
specification of the conditional covariance matrix $\Sigma_j(\theta)$ are unknown. So, to approximate the null distribution of this test statistic
we considered the following conditional resampling scheme, given the data $X_1, \ldots, X_n$:
 \begin{itemize} \itemsep=0pt
\item[(i)] Calculate $\widehat{\theta}_n=\widehat{\theta}_n(X_1, \ldots, X_n)$, the residuals  $\widetilde{\varepsilon}_1,\ldots,\widetilde{\varepsilon}_n$ and the test statistic $\widetilde{T}_{n,\gamma}^G=
\widehat{T}_{n,\gamma}^G(\widetilde{\varepsilon}_1,\ldots,\widetilde{\varepsilon}_n)$.

\item[(ii)] Generate vectors ${\varepsilon}_1^*,\ldots,{\varepsilon}_n^*$ i.i.d. from a N$_d(0,{ {\rm I}}_d)$ distribution. Let  $X_j^*=\Sigma_j^{1/2}(\widehat{\theta})\varepsilon_j^*$, $j=1,\ldots, n$.

\item[(iii)]    Calculate $\widehat{\theta}_n^*=\widehat{\theta}_n(X_1^*, \ldots, X_n^*)$, the residuals  $\widetilde{\varepsilon}_1^*,\ldots,\widetilde{\varepsilon}_n^*$, and approximate the null distribution of $\widetilde{T}^G_{n,\gamma}$ by means of the conditional distribution, given the data, of  $\widetilde{T}^{G*}_{n,\gamma}=\widetilde{T}^G_{n,\gamma}(\widetilde{\varepsilon}_1^*,\ldots,\widetilde{\varepsilon}_n^*)$.
\end{itemize}

In practice, the approximation in step (iii)   is carried out by generating a large number of bootstrap replications of the test statistic  $\widetilde{T}_{n,\gamma}^{G}$, for $b=1,\ldots, M$, whose empirical distribution function is used to estimate the null distribution of $\widetilde{T}^G_{n,\gamma}$.

In our simulation study, for the distribution of the innovations, we took  ${\varepsilon}_1, \ldots, {\varepsilon}_n$ i.i.d. from the distribution of ${\varepsilon}$ with
${\varepsilon}$ having a ($N$)  multivariate normal distribution, in order to study the level of the resulting test; and to study the power we considered the following heavy-tailed distributions:  $T_\theta$ and $GN_{\theta}$ (see \S8.1),  and the asymmetric exponential power distribution ($AEP$), whereby $(Z_1, \ldots, Z_d)^\top$, with $Z_1,\ldots,  Z_d$ i.i.d. from a univariate $AEP$ distribution  (\cite{zz:09})
 with parameters $\alpha=0.4$,  $p_1=1.182$ and $p_2=1.820$; these settings gave useful results in practical applications for the  errors in GARCH
type models. As in the previous subsection, we also calculated the HW-test for several values of $\beta$ with $Y_j$ replaced by $\widetilde{\varepsilon}_{j}$.

The parameters in the CCC-GARCH models were estimated by QMLE using the package {\tt ccgarch} of the language R.
Table \ref{tablagarch1}  reports the percentages of rejections for nominal significance level $\alpha=0.05$ and sample size $n=300$, for $r=0$ and $r=0.3$.
In order to reduce the computational burden we adopted the warp-speed method of  \cite{gpw:13} for evaluating the above resampling scheme.
With the warp-speed method, rather than computing critical points for each Monte Carlo sample, one resample is generated for each Monte Carlo sample,
and the resampling test statistic is computed for that sample. Then the resampling critical values for $\widetilde{T}^G_{n,\gamma}$ are computed from the empirical distribution
determined by the resampling replications of $\widetilde{T}_{n,\gamma}^{G*}$. In our simulations we took $10,000$ Monte Carlo samples for the level and $2,000$ for the power.
Looking at Table  \ref{tablagarch1}, we conclude that  the proposed test in most cases outperforms the HW-test by a wide margin. 

\begin{table} 
\centering
\caption{Percentage of rejection for nominal level $\alpha=0.05$ and $n=300$.} \label{tablagarch1}
\footnotesize
\begin{tabular}{|ccc|rrrrr|rrrrr|}
\hline
 &   &  &  \multicolumn{5}{c|}{Test based on $\widetilde{T}^G_{n,\gamma}$} & \multicolumn{5}{c|}{HW-test }\\ \cline{4-13}
 &   &  &   \multicolumn{5}{c|}{$\gamma$} & \multicolumn{5}{c|}{$\beta$}\\  \hline
    & $d$  & $r$   &  1.2 &  1.3&  1.4 &  1.5 & 2.0  & 0.5 & 1.0 & 1.5 & 2.0 & 2.5 \\ \hline
$N$ &  2     & 0.0 & 4.86  & 4.86  &  4.90 & 4.82 &  4.68 & 5.33 & 5.27 & 5.04 & 5.02 & 4.88 \vspace{-3pt}\\
    &        & 0.3 & 4.63  & 4.56  &  4.59 & 4.60 &  4.64 & 4.57 & 4.69 & 4.99 & 4.89 & 5.06 \vspace{-3pt}\\ 
    &  3     & 0.0 & 3.81  & 3.77  &  3.75 & 3.81 &  3.88 & 4.68 & 4.68 & 5.01 & 4.98 & 5.03 \vspace{-3pt}\\
    &        & 0.3 & 6.84  & 7.08  &  7.15 & 6.97 &  7.08 & 4.85 & 5.02 & 5.65 & 5.16 & 5.14 \vspace{-3pt}\\ \hline
$T_{10}$ & 2 & 0.0 & 89.15 & 89.05 & 88.85 & 88.70 & 89.05 & 4.40 & 17.00 & 27.75 & 29.10 &  27.70  \vspace{-3pt}\\
         &   & 0.3 & 81.30 & 81.15 & 81.10 & 81.40 & 81.95 & 5.10 & 18.75 & 26.50 & 27.25 & 26.00  \vspace{-3pt}\\
         & 3 & 0.0 & 95.80 & 95.90 & 95.80 & 95.80 & 96.05 & 6.15 & 35.10 & 46.90 & 43.20 & 37.05  \vspace{-3pt}\\
         &   & 0.3 & 16.10 & 13.15 & 12.90 & 13.25 & 15.35 & 9.10 & 37.05 & 45.85 & 43.75 & 37.25  \vspace{-3pt}\\ \hline
$GN_{1.65}$ & 2  & 0.0 & 48.70 & 48.60 & 49.00 & 48.85 & 48.60 &  5.50 &  7.20 & 11.20 & 13.00 & 13.5   \vspace{-3pt}\\
            &    & 0.3 & 48.15 & 47.80 & 47.30 & 47.20 & 46.25 &  4.75 &  7.40 & 12.10 & 13.80 & 14.25   \vspace{-3pt}\\
            & 3  & 0.0 & 59.10 & 59.80 & 59.35 & 59.20 & 58.00 &  4.10 &  7.55 & 11.65 & 12.15 & 12.85  \vspace{-3pt}\\
            &    & 0.3 & 54.55 & 54.00 & 54.25 & 54.45 & 53.80 &  4.35 &  8.15 & 11.85 & 12.70 & 11.95  \vspace{-3pt}\\  \hline
$AEP$ &  2   & 0.0 & 84.25 & 83.90 & 83.60 & 83.55 & 83.40 &  6.90 & 31.10 & 47.40 & 48.95 & 46.10  \vspace{-3pt}\\
      &      & 0.3 & 77.55 & 77.50 & 77.35 & 77.80 & 78.10 &  6.65 & 29.00 & 44.95 & 47.45 & 46.00 \vspace{-3pt}\\
      &  3   & 0.0 & 86.45 & 86.15 & 86.40 & 86.30 & 86.15 &  6.00 & 29.90 & 44.35 & 44.45 & 37.3  \vspace{-3pt}\\
      &      & 0.3 & 14.00 & 11.85 & 12.25 & 12.25 & 13.15 &  7.00 & 32.80 & 44.00 & 44.90 & 38.00          \vspace{-3pt}\\ \hline
\end{tabular}
\end{table}

\section{Conclusion}\label{sec_conclusion}
We prove new characterizations of the multivariate normal distribution. Based on these characterizations we suggest tests for normality with i.i.d.
data as well as tests for the null hypothesis of a classical Gaussian multivariate GARCH model. The new test statistics are simple to implement
and in limiting cases yield well-known moment-based criteria of goodness-of-fit.
Their asymptotic null distribution is studied and consistency against fixed alternatives is proved.
In a series of Monte Carlo experiments the new tests are implemented either based on asymptotic theory or by means of resampling schemes
showing that the resulting power is often higher than competitors.
For the GARCH case in particular, a  conditional resampling version of the procedure shows good power performance against certain heavy-tailed alternatives to the multivariate normal distribution.
%
%
%
%
%
\section{Appendix: Technical proofs}\label{sec_app}
Subsection \ref{Proofs} sketches the proofs of theoretical results. Subsection \ref{Auxiliary} contains some auxiliary results used in the proofs of our main theorems.
\subsection{Proofs} \label{Proofs}

\noindent {\sc Proof} of Theorem \ref{mainthm}.
Because of affine invariance of $T_n$, we
assume without loss of generality that $\BE X_1 = 0$ and $\BE[X_1 X_1^\top] = \textrm{I}_d$.
The main idea for showing convergence of the process $U_n(\cdot)= \sqrt{n}(M_n(\cdot)R_n(\cdot)-1)$ is as follows: Putting
\begin{eqnarray*}
A_n(t) & = & \exp\left(-\|t\|^2/2\right) \sqrt{n}\big{(}M_n(t) - \exp\left(\|t\|^2/2\right)\big{)},\\ \label{defbn}
B_n(t) & = & \exp\left(\|t\|^2/2\right) \sqrt{n}\big{(}R_n(t) - \exp\left(-\|t\|^2/2\right)\big{)},
\end{eqnarray*}
where $M_n$ and $R_n$ are given in (\ref{EMGF}) and (\ref{ECT}), respectively, the crucial observation, obtained after
straightforward computations, is the representation
\begin{equation}\label{repun}
U_n(t) = \frac{1}{\sqrt{n}} A_n(t)B_n(t) + A_n(t) + B_n(t), \quad t \in \RR^d.
\end{equation}
Hence, the program is to prove $A_n+B_n\vertk W$ in $\LL$. Since the first term on the right-hand side of $(\ref{repun})$
will turn out to be asymptotically negligible, Slutzky's lemma then gives the convergence $U_n \vertk W$.

The main problem of dealing with $A_n+B_n$ is that both $A_n$ and $B_n$ are sums of functions of the scaled residuals $Y_{n,1},\ldots,Y_{n,n}$
and not of the i.i.d. random vectors $X_1,\ldots,X_n$. To make the reasoning transparent, put
\begin{eqnarray*}
A_n^\circ (t) & = & \exp\left(-\|t\|^2/2\right) \sqrt{n}\big{(}M_n^\circ (t) - \exp\left(\|t\|^2/2\right)\big{)},\\
B_n^\circ (t) & = & \exp\left(\|t\|^2/2\right) \sqrt{n}\big{(}R_n^\circ (t) - \exp\left(-\|t\|^2/2\right)\big{)},
\end{eqnarray*}
where, in contrast to $R_n$ and $M_n$,
\begin{equation}\label{procrn0}
R_n^\circ (t) \ := \ \frac{1}{n} \sum_{j=1}^n \cos\left(t^\top X_j\right), \qquad M_n^\circ (t) \ := \ \frac{1}{n} \sum_{j=1}^n \exp\left(t^\top X_j\right), \quad t \in \RR^d,
\end{equation}
are based on $X_1,\ldots,X_n$. Notice that
$
A_n^\circ + B_n^\circ = n^{-1/2} \sum_{j=1}^n Z_j,
$
where
$
Z_j(t) \ = \ \e^{-\|t\|^2/2} \, \e^{t^\top X_j} - 1 + \e^{\|t\|^2/2}\cos\left(t^\top X_j\right)-1, t \in \RR^d.
$
Since $\BE[Z_1(t)] = 0$, $t \in \RR^d$, and (due to $\gamma >1$)  $\BE \|Z_1\|^2_{\textrm{L}^2} = \int_{\RR^d} \BE\left[ Z_1^2(t)\right] \, w_\gamma(t) \, \textrm{d}t < \infty$,
a Central Limit Theorem (CLT) in Hilbert spaces (see e.g., \cite{kmc:00}) shows that there is a centered Gaussian random element $W^\circ$ (say) of $\LL$, so that $A_n^\circ + B_n^\circ \vertk W^\circ.$

The idea now is to approximate each of the differences $A_n(t)-A_n^\circ(t)$ and $B_n(t)-B_n^\circ(t)$ by sums of i.i.d. random variables.
To this end, put
\begin{equation}\label{delta1}
\Delta_{n,j} = Y_{n,j} - X_j = \left(S_n^{-1/2}- \textrm{I}_d\right)X_j - S_n^{-1/2}\overline{X}_n, \qquad j=1,\ldots,n.
\end{equation}
A Taylor expansion gives
\begin{equation}\label{taylexp}
\e^{t^\top Y_{n,j}} - \e^{t^\top X_j} =
\e^{t^\top X_j} \! \left( \! t^\top \Delta_{n,j} + \frac{1}{2} \left(t^\top \Delta_{n,j}\right)^2 \exp\left(\Theta_{n,j} t^\top \Delta_{n,j} \right)\! \right),
\end{equation}
where $|\Theta_{n,j}| \le 1$. It follows that
\begin{equation}
A_n(t) - A_n^\circ(t) \ = \ \exp\left(-\|t\|^2/2\right) \, \sqrt{n} \, \big{(}M_n(t) - M_n^\circ(t)\big{)} \ = \   V_{n,1}(t) + V_{n,2}(t),
\end{equation}
where
\begin{eqnarray} \label{vn1}
V_{n,1}(t) & = &  \exp\left(-\|t\|^2/2\right) \, \frac{1}{\sqrt{n}} \sum_{j=1}^n \e^{t^\top X_j} t^\top \Delta_{n,j},\\ \label{vn2}
V_{n,2}(t) & = & \exp\left(-\|t\|^2/2\right) \, \frac{1}{2} \, \frac{1}{\sqrt{n}} \sum_{j=1}^n \e^{t^\top X_j} \left(t^\top \Delta_{n,j}\right)^2 \exp\left({\Theta_{n,j} t^\top \Delta_{n,j}}\right).
\end{eqnarray}
By some tedious and delicate estimations, we have $\|V_{n,2}\|^2_{\LL} =  o_\PP(1)$ (see Proposition \ref{asynegvn} in the Appendix),
and Proposition \ref{asyvn1} yields
\[
V_{n,1}(t) =  - \frac{1}{2} \, \frac{1}{\sqrt{n}} \sum_{j=1}^n t^\top \left(X_jX_j^\top - \textrm{I}_d\right)t - \frac{1}{\sqrt{n}} \sum_{j=1}^n t^\top X_j + o_\PP(1),
\]
where $o_\PP(1)$ refers to convergence in $\LL$.
Likewise, we have
\[
B_n(t) - B_n^\circ(t) \ = \ \exp\left(\|t\|^2/2\right) \, \sqrt{n} \, \left(R_n(t) - R_n^\circ(t)\right) \ = \   W_{n,1}(t) + W_{n,2}(t),
\]
where
\begin{eqnarray} \label{wn1}
W_{n,1}(t) & = &  - \, \exp\left(\|t\|^2/2\right) \, \frac{1}{\sqrt{n}} \sum_{j=1}^n \sin\left(t^\top X_j\right) t^\top \Delta_{n,j},\\ \nonumber
W_{n,2}(t) & = & \exp\left(\|t\|^2/2\right) \, \frac{1}{2} \, \frac{1}{\sqrt{n}} \sum_{j=1}^n \Psi_{n,j}  \left(t^\top \Delta_{n,j}\right)^2
\end{eqnarray}
and $|\Psi_{n,j}| \le 1$. Since
\[
|W_{n,2}(t)| \le \exp\left(\|t\|^2/2\right) \, \frac{\|t\|^2}{2} \,  \frac{1}{\sqrt{n}} \sum_{j=1}^n \|\Delta_{n,j}\|^2
\]
and $n^{-1/2} \sum_{j=1}^n \|\Delta_{n,j}\|^2 = o_\PP(1)$ (see \cite{hw:97}, p. 9), it follows that $\|W_{n,2}\|^2_{\LL} = o_\PP(1)$.
Moreover, Proposition \ref{asywn1} gives
\[
W_{n,1}(t) =   \frac{1}{2} \, \frac{1}{\sqrt{n}} \sum_{j=1}^n t^\top \left( X_jX_j^\top - \textrm{I}_d \right) t + o_\PP(1).
\]
Summarizing, we have
$
A_n(t) + B_n(t)  =  n^{-1/2} \sum_{j=1}^n Z_j(t) + o_\PP(1)$,
where
\[
Z_j(t)  =  \exp\left(- \|t\|^2/2\right) \e^{t^\top X_1} - 1 + \exp\left(\|t\|^2/2\right) \cos\left(t^\top X_j\right) - 1 - t^\top X_j.
\]
Since $\BE Z_j(t) =0$, $t\in \RR^d$, and, in view of the condition $\gamma >1$, $\BE \|Z_1\|^2_{\LL} < \infty$, a Hilbert space CLT (see \cite{kmc:00}) yields
$A_n + B_n \vertk W$ for a centered Gaussian element $W$ of $\LL$ having covariance kernel $C(s,t) = \BE[Z_1(s)Z_1(t)]$. Using
\begin{eqnarray*}
\BE\left[ \cos\left(s^\top X\right) \e^{t^\top X} \right] & = & \e^{(\|t\|^2 - \|s\|^2)/2} \, \cos \left(s^\top t\right), \\
\BE \left[ s^\top X \,  \e^{t^\top X}\right] & = & s^\top t \, \e^{\|t\|^2/2},\\
\BE \left[ \cos \left(s^\top X \right)  \cos\left( t^\top X\right) \right] & = & \frac{1}{2} \e^{- ( \|s\|^2 + \|t\|^2)/2} \, \left(\e^{s^\top t} + \e^{-s^\top t} \right) ,
\end{eqnarray*}
straightforward computations
show that $C(s,t)$ takes the form given in Theorem \ref{mainthm}. As a side-product of the derivations, we see that both $A_n$ and $B_n$ converge in distribution
and thus are tight sequences in $\LL$. Hence, the first term in (\ref{repun}) is $o_\PP(1)$, which completes the proof of Theorem \ref{mainthm}. \bewend

\medskip

\noindent {\sc Proof of (\ref{mean})}: Putting $I^{(0)}_\gamma:=I_\gamma(0,0)$, we have
\[
\mathbb E \|W\|^2_{\LL}=\frac{3}{2} I^{(0)}_{\gamma-1}+\frac{1}{2} I^{(0)}_{\gamma+1}-4 I^{(0)}_\gamma-J_{\gamma}+2 K_\gamma, \]
where $
J_\gamma=\int_{\mathbb R^d} \|t\|^2 \exp\left(-\gamma \|t\|^2\right) \textrm{d}t, \
K_\gamma=\int_{\mathbb R^d} \cos(\|t\|^2) \exp\left(-\gamma \|t\|^2\right) \textrm{d}t$.
Straightforward algebra gives $J_\gamma=d \: j^{(2)}_{\gamma} \left(j^{(0)}_{\gamma}\right)^{d-1}$,
where
\[
j^{(m)}_{\gamma}=\int_{-\infty}^\infty t^m \exp\left(-\gamma t^2\right) \textrm{d}t, \ m=0,2.
\]
Invoking
\[ \cos\left(\sum_{\ell=1}^d
\theta_\ell\right)=\sum_{q=0}^{\lfloor \frac{d}{2}\rfloor } (-1)^q \sum_{1\leq j_1
< j_2 ...< j_{2q} \leq d}  \ \prod_{k=1}^{2q} \sin \theta_{j_k} \
\prod_{\ell \neq j_1,j_2,...,j_{2q}}^{d} \cos \theta_\ell,\]
we obtain
\[
K_\gamma= \sum_{q=0}^{\lfloor \frac{d}{2} \rfloor } (-1)^q \: {d \choose 2q} \left(\kappa^{(s)}_{\gamma}\right)^{2q}\left(\kappa^{(c)}_{\gamma}\right)^{d-2q},
\]
where
$
\kappa^{(s)}_{\gamma}=\int_{-\infty}^\infty \sin(t^2) \exp\left(-\gamma t^2\right) \textrm{d}t, \ \ \kappa^{(c)}_{\gamma}=\int_{-\infty}^\infty \cos(t^2) \exp\left(-\gamma t^2\right) \textrm{d}t.
$
Since
\begin{eqnarray*}
j^{(0)}_{\gamma} & = & \sqrt{\frac{\pi}{\gamma}}, \quad  j^{(2)}_{\gamma}=\frac{1}{2\gamma} \sqrt{\frac{\pi}{\gamma}},\\
\kappa^{(s)}_{\gamma} & = & \sqrt{\frac{\pi}{\gamma}}\frac{\sin\left((1/2)\tan^{-1}(1/\gamma)\right)}{(1+(1/\gamma^2))^{1/4}}, \quad \ \kappa^{(c)}_{\gamma}=\sqrt{\frac{\pi}{\gamma}}\frac{\cos\left((1/2)\tan^{-1}(1/\gamma)\right)}{(1+(1/\gamma^2))^{1/4}},
\end{eqnarray*}
the result follows by simple algebra. \bewend

\medskip

\noindent {\sc Proof} of Theorem \ref{mainthmG}.
Notice that $U_n^G$ satisfies an equality similar to that in (\ref{repun}) with $A_n$ and $B_n$ replaced with $A_n^G$ and $B_n^G$, respectively, where
\[A_n^G(t)=\exp(-\|t\|^2/2)\frac{1}{\sqrt{n}}\sum_{j=1}^n \left[\exp\left(t^\top \widetilde{\varepsilon}_j(\widehat{\theta}_n)\right)-\exp(\|t\|^2/2)\right], \]
\[B_n^G(t)=\exp(\|t\|^2/2)\frac{1}{\sqrt{n}}\sum_{j=1}^n \left[\cos\left(t^\top \widetilde{\varepsilon}_j(\widehat{\theta}_n)\right)-\exp(-\|t\|^2/2)\right]. \]
To prove the result we will demonstrate that
\begin{equation} \label{estoes}
\begin{array}{rcl}
A_n^G(t)+B_n^G(t)& = & \frac{1}{\sqrt{n}}\sum_{j=1}^n \! \left[\! \exp \! \left( \! -\frac{\|t\|^2}{2}\right)\exp(t^\top \varepsilon_j)+\exp \! \left(\frac{\|t\|^2}{2}\right)\cos(t^\top \varepsilon_j)\! - \! 2\right]\\ & & +r_n(t),
 \end{array}
\end{equation}
with $\|r_{n}\|_{\textrm{L}^2}=o_{\mathbb{P}}(1)$,  and the result will follow from the CLT in Hilbert spaces.
With this aim, we first introduce some notation. Let $\varepsilon_j(\theta)=\Sigma_j^{-1/2}(\theta)X_j$. Notice that  $\varepsilon_j(\theta_0)=\varepsilon_j$. Let
\[
W_j=\varepsilon_j+\sum_{k=1}^v\left.\frac{\partial}{\partial \theta_k} \varepsilon_j(\theta)\right|_{\theta=\theta_0}(\widehat{\theta}_{nk}-\theta_{0k}),
\]
where $\widehat{\theta}_n=(\widehat{\theta}_{n1}, \ldots, \widehat{\theta}_{nv})^\top$, $\theta_0=(\theta_{01}, \ldots, \theta_{0v})^\top$, and $\Delta_{n,j}=\varepsilon_j(\widehat{\theta}_n)-W_j$.
 Then
\[
\widetilde{\varepsilon}_j(\widehat{\theta}_n)-W_j=\widetilde{\varepsilon}_j(\widehat{\theta}_n)-{\varepsilon}_j(\widehat{\theta}_n)+\Delta_{n,j}.
\]
Let $A_{jk}(\theta)=\Sigma_j^{-1/2}(\theta)\frac{\partial}{\partial \theta_k} \Sigma_j^{1/2}(\theta)$ and $\mu_k=\mathbb{E}[A_{jk}(\theta)]$.
To show  (\ref{estoes}) we prove that
\begin{itemize} \itemsep=0pt
\item[(a.1)] $\exp(-\|t\|^2/2)\frac{1}{\sqrt{n}}\sum_{j=1}^n\exp(t^\top W_j)=\frac{1}{\sqrt{n}}\sum_{j=1}^nV_j(t)+r_{n,1}(t)$, with  $$V_j(t)=\exp(-\|t\|^2/2+t^\top\varepsilon_j)-\sum_{k=1}^vt^\top\mu_{k}t\sqrt{n}(\theta_{nk}-\theta_{0k}), \quad
\|r_{n,1}\|_{\textrm{L}^2}=o_{\mathbb{P}}(1),$$

\item[(a.2)] $\exp(\|t\|^2/2)\frac{1}{\sqrt{n}}\sum_{j=1}^n\cos(t^\top W_j)=\frac{1}{\sqrt{n}}\sum_{j=1}^nZ_j(t)+r_{n,2}(t)$, with  $$Z_j(t)=\exp(\|t\|^2/2)\cos(t^\top\varepsilon_j)+\sum_{k=1}^vt^\top\mu_{k}t\sqrt{n}(\theta_{nk}-\theta_{0k}), \quad
\|r_{n,2}\|_{\textrm{L}^2}=o_{\mathbb{P}}(1),$$

\item[(b.1)] $\|r_{n,3}\|_{\textrm{L}^2}=o_{\mathbb{P}}(1)$, where $r_{n,3}(t)=\textrm{e}^{-\|t\|^2/2}\frac{1}{\sqrt{n}}\sum_{j=1}^n\exp(t^\top W_j)\left(\exp(t^\top \Delta_{n,j})-1\right),$

\item[(b.2)] $\|r_{n,4}\|_{\textrm{L}^2}=o_{\mathbb{P}}(1)$, where $r_{n,4}(t)=\textrm{e}^{\|t\|^2/2}\frac{1}{\sqrt{n}}\sum_{j=1}^n\cos(t^\top W_j)\left(\exp(t^\top \Delta_{n,j})-1\right),$

\item[(c.1)] $\|r_{n,5}\|_{\textrm{L}^2}=o_{\mathbb{P}}(1)$, where $r_{n,5}(t)=\textrm{e}^{-\|t\|^2/2}\frac{1}{\sqrt{n}}\sum_{j=1}^n\left[\exp(\widetilde{\varepsilon}_j(\widehat{\theta}_n))-\exp({\varepsilon}_j(\widehat{\theta}_n))\right],$

\item[(c.2)] $\|r_{n,6}\|_{\textrm{L}^2}=o_{\mathbb{P}}(1)$, where $r_{n,6}(t)=\textrm{e}^{\|t\|^2/2}\frac{1}{\sqrt{n}}\sum_{j=1}^n\left[\cos(\widetilde{\varepsilon}_j(\widehat{\theta}_n))-\cos({\varepsilon}_j(\widehat{\theta}_n))\right].$
\end{itemize}

\noindent \underline{Proof of (a.1)}. Notice that $\left.\frac{\partial}{\partial \theta_k} \varepsilon_j(\theta)\right|_{\theta=\theta_0}=-A_{jk}(\theta_0)\varepsilon_j$. A Taylor expansion gives
\begin{eqnarray*}
\frac{1}{\sqrt{n}}\sum_{j=1}^n\exp(t^\top W_j) & = & \frac{1}{\sqrt{n}}\sum_{j=1}^n \exp(t^\top \varepsilon_j)\\
 & & -\frac{1}{\sqrt{n}}\sum_{j=1}^n \exp(t^\top \varepsilon_j)
t^\top\sum_{k=1}^vA_{jk}(\theta_0)(\widehat{\theta}_{nk}-\theta_{0k})\varepsilon_j+ \frac{R_{n,1}(t)}{2},
\end{eqnarray*}
with $|R_{n,1}|\leq n^{-1/2} R_{n,1,1} R_{n,1,2}$, where
\begin{eqnarray*}
R_{n,1,1}(t) & = & \exp\left(\|t\| \max_j\|\sum_{k=1}^vA_{jk}(\theta_0)\| \|\widehat{\theta}_n-\theta_0\| \max_j\|\varepsilon_j\|\right),\\
R_{n,1,2}(t) & = & \sum_{j=1}^n \exp(t^\top \varepsilon_j)\left(t^\top\sum_{k=1}^vA_{jk}(\theta_0)(\widehat{\theta}_{nk}-\theta_{0k})\varepsilon_j\right)^2.
\end{eqnarray*}
We have
\[
\exp(-\frac{\|t\|^2}{2})\frac{1}{\sqrt{n}}\sum_{j=1}^n \textrm{e}^{t^\top \varepsilon_j}
t^\top\sum_{k=1}^vA_{jk}(\theta_0)(\widehat{\theta}_{nk} \!- \! \theta_{0k})\varepsilon_j=t^\top\sum_{k=1}^v\mu_{k}\sqrt{n}(\theta_{nk}\! -\! \theta_{0k})t\]\[+R_{n,1,3}(t),
\]
where
\begin{eqnarray*}
R_{n,1,3}(t) & = & \sum_{k=1}^v \sqrt{n}(\widehat{\theta}_{nk} -\theta_{0k})R_{n,1,3,k}(t), \quad
R_{n,1,3,k}(t)=\frac{1}{n}\sum_{j=1}^nV_{kj}(t),\\
V_{kj}(t) & = & \exp(-\|t\|^2/2)\exp(t^\top \varepsilon_j)t^\top A_{jk}(\theta_0)\varepsilon_j-t^\top\mu_{k}t.
\end{eqnarray*}
By noting that  $\mathbb{E}[V_{kj}]=0$, $\mathbb{E}[\langle V_{kj}, V_{kr}\rangle]=0$, $\forall j \neq r$, $\forall k$, and $\mathbb{E}\|V_{kj}\|^2_{\textrm{L}^2}<\infty$, it follows that $\|R_{n,1,3,k}\|_{\textrm{L}^2}=o_{\mathbb{P}}(1)$, $1\leq k \leq v$, and thus $\|R_{n,1,3}\|_{\textrm{L}^2}=o_{\mathbb{P}}(1)$.
Similar calculations show $\|\exp(-\|t\|^2/2)R_{n,1,2}\|_{\textrm{L}^2}=O_{\mathbb{P}}(1)$.
Let $\alpha>0$ so that $1/q<\alpha<1/2$, where $q$ is as in (A.6). Then
\begin{eqnarray*}
\mathbb{P}\left(\max_j \Big{\|}\sum_{k=1}^v A_{jk}(\theta_0)\Big{\|}> n^{\alpha}\right) & \leq  & n \, \mathbb{P}\left(\Big{\|}\sum_{k=1}^v A_{jk}(\theta_0)\Big{\|}> n^{\alpha}\right) \\
& \leq  & n \, \frac{\mathbb{E}(\|\sum_{k=1}^v A_{jk}(\theta_0)\|^r)}{n^{r\alpha}}\to 0 \,\,\,\, \textrm{as} \,\,\,\, n \to \infty.
\end{eqnarray*}
This convergence together with (A.1) and Proposition \ref{propo1} imply
 \[\max_j\Big{\|}\sum_{k=1}^vA_{jk}(\theta_0)\Big{\|} \|\widehat{\theta}_n-\theta_0\| \max_j\|\varepsilon_j\|=o_{\mathbb{P}}\left(\frac{\sqrt{\log(n)}}{n^{0.5-\alpha}}\right).\]
 Therefore, Proposition \ref{propint} gives $\|\exp(-\|t\|^2/2)R_{n,1,1}\|_{\textrm{L}^2}=O_{\mathbb{P}}(1)$. As a consequence, $\|\exp(-\|t\|^2/2)R_{n,1}\|_{\textrm{L}^2}=o_{\mathbb{P}}(1)$, and  (a.1) follows.

 \noindent \underline{Proof of (a.2)}. The proof is similar to that of (a.1), so we omit it.

 \noindent \underline{Proof of (b.1)}. Observe that
 $\exp(t^\top \Delta_{n,j})-1=t^\top \Delta_{n,j}\exp(\alpha_{n,j}t^\top \Delta_{n,j})$
 for some $\alpha_{n,j}\in (0,1)$, and
 \[\Delta_{n,j}=\sum_{k,l=1}^v\left. \frac{\partial^2}{\partial \theta_k \partial \theta_l} \varepsilon_j(\theta) \right|_{\theta=\widetilde{\theta}_n}(\widehat{\theta}_{nk}-\theta_{0k})(\widehat{\theta}_{nl}-\theta_{0l})\]
 for some $\widetilde{\theta}_n$ between $\widehat{\theta}_{n}$ and $\theta_0$. Now (A.1) and (A.6) yield
 $\Delta_{n,j}\| \leq D_j \|\varepsilon_j\|\| \widehat{\theta}_{n}-\theta_0\|^2$ for large enough $n$, where
 $ \mathbb{E} D_j^2<\infty$. The Cauchy--Schwarz inequality gives
 \[|r_{n,3}(t)|\leq r_{n,3,1}(t)^{1/2} r_{n,3,2}(t)^{1/2},\]
 where
 \[ r_{n,3,1}(t)=\textrm{e}^{-\|t\|^2/2} \frac{1}{\sqrt{n}}\sum_{j=1}^n \exp(t^\top W_j), \]
 \[ r_{n,3,2}(t)=\textrm{e}^{-\|t\|^2/2} \|t\|^2\|\sqrt{n}(\widehat{\theta}_n-\theta_0)\|^4\frac{1}{\sqrt{n}n^2}\! \sum_{j=1}^n \! D_j^2\|\varepsilon_j\|^2\exp(2\|t\| \| D_j\| \|\varepsilon_j\|
 \| \widehat{\theta}_{n}\! -\! \theta_0\|^2).
 \]
Proceeding as in the proof of (a.1), it can be seen that $r_{n,3,1}=O_{\mathbb{P}}(1)$ (in $\textrm{L}^2$) and that $r_{n,3,2}=o_{\mathbb{P}}(1)$ (in $\textrm{L}^2$), showing that
$r_{n,3}=o_{\mathbb{P}}(1)$ (in $\textrm{L}^2$).

\noindent \underline{Proof of (b.2)}. The proof is similar to that of (b.1) and is thus omitted.

\noindent \underline{Proof of (c.1)}. Let $\Lambda_{n,j}=\widetilde{\varepsilon}_j(\widehat{\theta}_n)-{\varepsilon}_j(\widehat{\theta}_n)=
\widetilde{\Sigma}_j^{-1/2}(\widehat{\theta}_n) \left(\Sigma_j^{1/2}(\widehat{\theta}_n)-\widetilde{\Sigma}^{1/2}_j(\widehat{\theta}_n)\right) \Sigma_j^{-1/2}(\widehat{\theta}_n)X_j$.
A Taylor expansion yields
\begin{equation} \label{aux2}
\exp\left(t^\top\widetilde{\varepsilon}_j(\widehat{\theta}_n)\right) -\exp\left(t^\top{\varepsilon}_j(\widehat{\theta}_n)\right)=t^\top\Lambda_{n,j}\exp\left(t^\top{\varepsilon}_j(\widehat{\theta}_n)+\alpha_{n,j}t^\top\Lambda_{n,j}\right)
\end{equation}
for some  $\alpha_{n,j} \in (0,1)$.  By the Cauchy--Schwarz inequality  and (\ref{aux2}),
\[\exp(-\|t\|^2/2)\left|\frac{1}{\sqrt{n}}\sum_{j=1}^n \left[\exp\left(t^\top\widetilde{\varepsilon}_j(\widehat{\theta}_n)\right) -\exp\left(t^\top{\varepsilon}_j(\widehat{\theta}_n)\right)\right] \right|\leq R_{n,1}(t)^{1/2} R_{n,2}(t)^{1/2},
\]
where
\[R_{n,1}(t)=\exp(-\|t\|^2/2)\frac{1}{\sqrt{n}}\sum_{j=1}^n \exp\left(2t^\top{\varepsilon}_j(\widehat{\theta}_n)\right),\]
\[R_{n,2}(t)=\exp(-\|t\|^2/2)\frac{1}{\sqrt{n}}\sum_{j=1}^n \|t\|^2 \|\Lambda_{n,j} \|^2 \exp(2\|t\| \|\Lambda_{n,j} \|).\]
Proceeding as in (a.1), we obtain $\|R_{n,1}\|_2 =O_\mathbb{P}(1)$. We next show $\|R_{n,2}\|_{\textrm{L}^2}=o_{\mathbb{P}}(1)$, which yields (c.1). To this end,
put $\Lambda_{n,j,k} =\|\Lambda_{n,j} \|+ \|\Lambda_{n,k} \|$. From Proposition \ref{propint},
\[\|R_{n,2}\|^2_{\textrm{L}^2}\leq  \frac{1}{n}\sum_{j,k=1}^n\|\Lambda_{n,j} \|^2 \|\Lambda_{n,k} \|^2\left(K_0+K_1\Lambda_{n,j,k} ^{d+3}\right)
\exp\left(\frac{\Lambda_{n,j,k} ^2}{\gamma+1}\right). \]
Taking into account that for any $x,y \in \mathbb{R}$ and $r >0$
\begin{equation} \label{aux0}
|x+y|^r=c_r(|x|^r+|y|^r),
\end{equation}
where $c_r =1$ if $0<r \le 1$ and $c_r=2^{r-1}$, otherwise, we have
\[\|R_{n,2}\|^2_{\textrm{L}^2}  \leq   C \left\{\frac{1}{\sqrt{n}}\sum_{j=1}^n\|\Lambda_{n,j} \|^2 \exp\left(\frac{2\|\Lambda_{n,j}\|^2}{\gamma+1}\right)\right\}^2\]
\[+C\left\{\frac{1}{\sqrt{n}}\sum_{j=1}^n\|\Lambda_{n,j} \|^2 \exp\left(\frac{2\|\Lambda_{n,j}\|^2}{\gamma+1}\right)\right\}
\left\{\frac{1}{\sqrt{n}}\sum_{j=1}^n\|\Lambda_{n,j} \|^{d+5} \exp\left(\frac{2\|\Lambda_{n,j}\|^2}{\gamma+1}\right)\right\}.
\]
We show that each expression within curly brackets is $o_{\mathbb{P}}(1)$.
From (A.2) and (A.3),
\begin{equation} \label{aux1}
\|\Lambda_{n,j}\| \leq C\rho^j\|X_j\|.
\end{equation}
As a consequence, (A.4) implies
\begin{equation} \label{aux3}
\Lambda_{n,j} \to 0 \quad \textrm{a.s.}
\end{equation}
as $j\to \infty$ for each $n$, because the upper bound in (\ref{aux1}) does not depend on $n$ (see Exercise 7.2 in \cite{fz:10}.
From (\ref{aux3}), it follows that for each $M>1$ and each $\omega \in \Omega$, there is an integer $j_0=j_0(\omega, M)$ so that $\exp\left(\frac{2\|\Lambda_{n,j}\|^2}{\gamma+1}\right) \leq M$
for each $n$ and each $j,k >j_0$. For $n> n_0$ and $r=2$ or $r=d+5$,
\[\frac{1}{\sqrt{n}}\sum_{j=1}^n\|\Lambda_{n,j} \|^r \exp\left(\frac{2\|\Lambda_{n,j}\|^2}{\gamma+1}\right)\leq T_1+T_2,
\]
where
\[T_1=\frac{n_0}{\sqrt{n}}\max_{1\leq j \leq n_0}\|\Lambda_{n,j} \|^r \exp\left(\frac{2\|\Lambda_{n,j}\|^2}{\gamma+1}\right), \quad
T_2=\frac{M}{\sqrt{n}}\sum_{j=n_0+1}^n\|\Lambda_{n,j} \|^r .\]
Let $0<\alpha<1/2$, and put $\log^+(x)=\max\{0,\log(x)\}$. Then
\[
\mathbb{P}\left(\max_{1\leq j \leq n_0}\|\Lambda_{n,j} \|^r \exp\left(\frac{2\|\Lambda_{n,j}\|^2}{\gamma+1}\right) \! > \! n^\alpha \! \right) \leq
\sum_{j=1}^{n_0}\mathbb{P} \! \left(\! \|\Lambda_{n,j} \|^r \exp\left(\frac{2\|\Lambda_{n,j}\|^2}{\gamma+1}\right)\! > \! n^\alpha\right)
\]
\[
\leq C \sum_{j=1}^{n_0}\frac{\mathbb{E}[\log^+(\|\Lambda_{n,j} \|)^{\zeta/2}] + \mathbb{E} \|\Lambda_{n,j} \|^\zeta}{(\log(n))^{\zeta/2}}
\leq C \frac{ \mathbb{E} \|X_1\|^{\zeta}\sum_{j=1}^\infty \rho^{j\zeta}}{(\log(n))^{\zeta/2}},
\]
which implies $\max_{1\leq j \leq n_0}\|\Lambda_{n,j} \|^r \exp\left(\frac{2\|\Lambda_{n,j}\|^2}{\gamma+1}\right)=o_{\mathbb{P}}(n^{\alpha})$ and thus $T_1=o_{\mathbb{P}}(1)$.
From (A.4) and (\ref{aux0}), for any $0<\zeta\leq \min \{1,\varsigma\}$ , we have
\[ \mathbb{E} \left(\sum_{j=n_0+1}^n\|\Lambda_{n,j} \|^r\right)^{\zeta/r} \leq \sum_{j=n_0+1}^n \mathbb{E} \|\Lambda_{n,j} \|^{\zeta} \leq \mathbb{E} \|X_1\|^{\zeta}\sum_{j=n_0+1}^\infty \rho^{j\zeta}<\infty. \]
Since $\sum_{j=n_0+1}^n\|\Lambda_{n,j} \|^r$ has finite moment of order $\zeta/r$, it is finite almost surely. Thus $T_2 \to 0$ a.s. as $n\to \infty$, which
completes the proof of   (c.1).

 \noindent \underline{Proof of (c.2)}. The proof is similar to that of (c.1) and is thus omitted. \bewend

\subsection{Some auxiliary results} \label{Auxiliary}
\begin{prop} \label{propo1}
Let $X,X_1,X_2, \ldots $ be i.i.d. $d$-variate random vectors having the normal distribution N$_d(0,\textrm{I}_d)$, and put
$F_n := \max_{1 \le j \le n} \|X_j\|$. Then
\[
\lim_{n \to \infty} \mathbb{P} \left(4 \sqrt{2 \log n}\left(F_n - a_n \right) \le t \right) \ = \ \exp(-\exp(-t)), \qquad t \in \mathbb{R},
\]
where
\[
a_n = \sqrt{2 \log n} + \frac{(d-2)\log \log n}{2 \sqrt{2 \log n}} - \frac{\log \Gamma(d/2)}{\sqrt{2 \log n} }.
\]
\end{prop}

\noindent {\sc Proof.}
Since $\|X\|^2$ has a $\chi^2_d$-distribution and thus a Gamma distribution, the distribution of $F_n^2$ is in the maximum domain of attraction of the Gumbel law. From p. 156 of  
\cite{ekm:97} we therefore obtain
$
2(F_n^2 - d_n) \ \vertk \ G,
$
where
\[
d_n \ = \ 2 \left( \log n + \left( \frac{d}{2} - 1 \right) \log \log n - \log \Gamma(d/2) \right)
\]
and $G$ has a Gumbel distribution.  Using $F_n^2- d_n = (F_n - \sqrt{d_n})(F_n + \sqrt{d_n})$ together with
$(F_n + \sqrt{d_n})/\sqrt{d_n} \stk 2$ and Sluzky's lemma, we have $4 \sqrt{d_n}(F_n - \sqrt{d_n}) \vertk Z$.
Upon noting that $\sqrt{d_n} = a_n + o(1/\sqrt{2 \log n})$ as $n \to \infty$, the assertion follows.
\bewend

\begin{prop} \label{propint}
Let $k \in \NN\cup \{0\}$ and $\alpha >0$. Then
\[
\int_{\RR^d} \|t\|^k \exp\left(-(1+\gamma)\|t\|^2 + 2 \, \alpha \, \|t\|\right) \, \textrm{d}t \ \le \ \left(K_0 + K_1 \, \alpha^{d+k-1}\right)
\exp\left(\frac{\alpha^2}{1+\gamma}\right)
\]
for some constants $K_0, K_1$ that depend only on $d$, $k$ and $\gamma$.
\end{prop}

\noindent {\sc Proof.} Using spherical coordinates, the integral equals
\[
C_d \, \exp\left(\frac{\alpha^2}{1+\gamma}\right) \int_0^\infty r^{k+d-1} \exp\left(-(1+\gamma)\left(r- \frac{\alpha}{1+\gamma}\right)^2\right) \, \textrm{d}r,
\]
where $C_d$ is a constant that depends only on $d$. This last integral, in turn, is equal to
$\frac{1}{2} \sigma \sqrt{2\pi} \BE\left[|N|^{k+d-1}\right]$, where $\sigma^2 = (2(1+\gamma))^{-1}$ and $N$ has the normal distribution N$(\mu,\sigma^2)$,
where $\mu = \alpha/(1+\gamma)$. From this, the result follows readily. \bewend

\begin{prop} \label{asynegvn}
We have
$
\|V_{n,2}\|^2_{\LL} \ = \ o_\PP(1),
$
where $V_{n,2}$ is given in (\ref{vn2}).
\end{prop}

\noindent {\sc Proof.} Putting
$F_n := \max_{1\le j \le n} \|X_j\|$, (\ref{delta1}) gives
\begin{equation}\label{tdeltaab}
|t^\top \Delta_{n,j}|  \le  \|t\| \, \Lambda_n,
\end{equation}
where
\begin{equation}\label{lambdan}
\Lambda_n \ = \ \|S_n^{-1/2}- \textrm{I}_d\|_2 \, F_n + \|S_n^{-1/2}\|_2 \, \|\overline{X}_n\|.
\end{equation}
On the other hand, we have
$
(t^\top \Delta_{n,j})^2 \le 2 (t^\top (S_n^{-1/2}- \textrm{I}_d) X_j )^2 + 2 (t^\top S_n^{-1/2}\overline{X}_n)^2,
$
which implies $
0 \le V_{n,2}(t) \le V_{n,2,1}(t) + V_{n,2,2}(t)$,
where
\begin{eqnarray}\nonumber
V_{n,2,1}(t) & = & \e^{-\|t\|^2/2}\, \frac{1}{\sqrt{n}} \sum_{j=1}^n \e^{t^\top X_j} \left(t^\top \left(S_n^{-1/2}- \textrm{I}_d\right) X_j \right)^2 \exp\left( \Theta_{n,j}t^\top \Delta_{n,j}\right),\\ \label{vn22}
V_{n,2,2}(t) & = & \e^{-\|t\|^2/2}\, \frac{1}{\sqrt{n}} \sum_{j=1}^n \e^{t^\top X_j} \left(t^\top S_n^{-1/2}\overline{X}_n\right)^2
\exp\left(\Theta_{n,j}t^\top \Delta_{n,j}\right).
\end{eqnarray}
Since $\|V_{n,2}\|^2_{\LL} \le 2 \|V_{n,2,1}\|^2_{\LL} + 2 \|V_{n,2,2}\|^2_{\LL}$, it suffices to prove that each of the last two summands is $o_\PP(1)$.
To tackle $V_{n,2,1}$, notice that
\[
\left(t^\top \left(S_n^{-1/2}- \textrm{I}_d\right) X_j \right)^2 \le \|t\|^2 \, \|S_n^{-1/2}-\textrm{I}_d \|^2_2 \, \|X_j\|^2
\]
and hence, invoking (\ref{tdeltaab}),
\[
V_{n,2,1}(t) \le \|S_n^{-1/2}-\textrm{I}_d\|^2_2 \, \|t\|^2 \, \e^{-\|t\|^2/2} \, \frac{1}{\sqrt{n}} \sum_{j=1}^n \e^{t^\top X_j} \, \|X_j\|^2 \, \e^{\|t\|\Lambda_n}.
\]
Consequently,
\[
\|V_{n,2,1}\|^2_{\LL} \le \|S_n^{-1/2}\! -\! \textrm{I}_d\|^4_2 \, \frac{1}{n} \sum_{i,j=1}^n \!  \|X_i\|^2\, \|X_j\|^2 \! \int_{\RR^d}
\|t\|^4 \e^{-(1+\gamma)\|t\|^2} \, \e^{t^\top (X_i+X_j)+2 \|t\| \Lambda_n} \, \textrm{d}t.
\]
Now, using $|t^\top (X_i+X_j)| \le 2 \|t\| \, F_n$ and Proposition \ref{propint} with $k=4$ and $\alpha = F_n + \Lambda_n$, it follows that
\begin{equation}\label{finalab}
\|V_{n,2,1}\|^2_{\LL} \le \|S_n^{-1/2}\! -\! \textrm{I}_d\|^4_2 \cdot n \cdot \frac{1}{n^2} \sum_{i,j=1}^n \|X_i\|^2\, \|X_j\|^2 \, (K_0+K_1\Gamma_n^{d+3})  \, \exp\left(\frac{\Gamma_n^2}{1+\gamma}\right),
\end{equation}
where $\Gamma _n = F_n + \Lambda_n$. From Proposition \ref{propo1}, we obtain
\[
F_n = \sqrt{2 \log n} + \frac{(d-2)\log \log n}{2 \sqrt{2 \log n}} + O_\PP\left(\frac{1}{\sqrt{\log n}}\right).
\]
In view of $\|S_n^{-1/2}-\textrm{I}_d\|_2 = O_\PP(n^{-1/2})$ and $\|\overline{X}_n\| = O_\PP(n^{-1/2})$, we see that
$\Lambda_n$ figuring in (\ref{lambdan}) is of order $O_\PP((\log n/n)^{1/2})$ and thus
\[
\Gamma_n = \sqrt{2 \log n} + \frac{(d-2)\log \log n}{2 \sqrt{2 \log n}} + O_\PP\left(\frac{1}{\sqrt{\log n}}\right).
\]
Hence, $\Gamma_n^2 =  2 \log n + (d-2)\log \log n + O_\PP(1)$ and therefore
\[
\frac{\Gamma_n^2}{1+ \gamma} = \log \left(n^{2/(1+\gamma)}\right) + \log \left(\left( \log n\right)^{(d-2)/(1+\gamma)} \right) + O_\PP(1).
\]
It follows that the rightmost factor figuring in (\ref{finalab}) is
\begin{equation} \label{darstellgamnq}
\exp\left(\frac{\Gamma_n^2}{1+\gamma}\right) = n^{2/(1+\gamma)} \, \left(\log n \right)^{(d-2)/(1+\gamma)} \, O_\PP(1).
\end{equation}
Since $\|S_n^{-1/2}-\textrm{I}_d\|^4_2 = O_\PP(n^{-2})$, $n^{-2} \sum_{i,j=1}^n \|X_i\|^2\, \|X_j\|^2= O_\PP(1)$ and  $\Gamma_n^{d+3}$ is of order
$O_\PP((\log n)^{(d+3)/2})$, (\ref{darstellgamnq}) yields
\[
\|V_{n,2,1}\|^2_{\LL} = O_\PP\left(n^{(1-\gamma)/(1+\gamma)}\right) \cdot \left( \log n\right)^{\frac{d-2}{1+\gamma} + \frac{d+3}{2}}
\]
and thus $\|V_{n,2,1}\|^2_{\LL} = o_\PP(1)$, since $\gamma >1$.

To show that $\|V_{n,2,2}\|^2_{\LL} =o_\PP(1)$, where $V_{n,2,2}$ is given in (\ref{vn22}), a similar reasoning as above
yields
$
\|V_{n,2,2}\|^2_{\LL} = O_\PP(n^{-1}) \cdot \int_{\RR^d} \|t\|^4 \e^{-(1+\gamma)\|t\|^2} \, \e^{2\|t\|(F_n + \Lambda_n)}\, \textrm{d}t$.
Hence, $\|V_{n,2,2}\|^2_{\LL}$ is of the same order as $\|V_{n,2,1}\|^2_{\LL}$, which completes the proof.
\bewend

\begin{prop} \label{asyvn1}
We have
\[
V_{n,1}(t) = - \frac{1}{2} \, \frac{1}{\sqrt{n}} \sum_{j=1}^n t^\top \left(X_jX_j^\top - \textrm{I}_d\right)t - \frac{1}{\sqrt{n}} \sum_{j=1}^n t^\top X_j + o_\PP(1),
\]
where $V_{n,1}$ is given in (\ref{vn1}) and $o_\PP(1)$ is with respect to $\LL$.
\end{prop}

\noindent {\sc Proof.} From the definition of $\Delta_{n,j}$, we have
$
V_{n,1}(t)  =  V_{n,1,1}(t)\! -\! V_{n,1,2}(t)\! -\!  V_{n,1,3}(t) ,
$
where
\begin{eqnarray*}
V_{n,1,1}(t) & = & e^{-\|t\|^2/2} \, \frac{1}{\sqrt{n}} \sum_{j=1}^n \e^{t^\top X_j} t^\top \left(S_n^{-1/2} - \textrm{I}_d\right) X_j,\\
V_{n,1,2}(t) & = & e^{-\|t\|^2/2} \, \frac{1}{\sqrt{n}} \sum_{j=1}^n \e^{t^\top X_j} t^\top \left(S_n^{-1/2} - \textrm{I}_d\right) \overline{X}_n,\\
V_{n,1,3}(t) & = & e^{-\|t\|^2/2} \, \frac{1}{\sqrt{n}} \sum_{j=1}^n \e^{t^\top X_j} t^\top \overline{X}_n.
\end{eqnarray*}
Since $|t^\top \left(S_n^{-1/2} - \textrm{I}_d\right) \overline{X}_n| \le \|t\|\, \|S_n^{-1/2}-\textrm{I}_d\|_2 \, \|\overline{X}_n\| = \|t\| \, O_\PP(n^{-1})$,
it is readily seen that $\|V_{n,1,2}\|^2_{\LL} = o_\PP(1)$. Now,
\[
V_{n,1,3}(t) = \e^{-\|t\|^2/2}  \left(M_n^\circ (t) - m(t)\right) \frac{1}{\sqrt{n}} \sum_{k=1}^n t^\top X_k + \frac{1}{\sqrt{n}} \sum_{j=1}^n t^\top X_j.
\]
with $M_n^\circ(t)$ given in (\ref{procrn0}) and $m(t) = \exp(\|t\|^2/2)$. Taking expectations and using Fubini's Theorem,
it follows that $V_{n,1,3}(t) = n^{-1/2} \sum_{j=1}^n t^\top X_j + o_\PP(1)$, where $o_\PP(1)$ is understood with respect to $\LL$.
Finally, we use the relation
\begin{equation}\label{snminushalb}
\sqrt{n} \left(S_n^{-1/2}- \textrm{I}_d \right) = - \frac{1}{2 \sqrt{n}} \sum_{k=1}^n \left(X_kX_k^\top - \textrm{I}_d \right) + O_\PP\left(n^{-1/2}\right)
\end{equation}
(see display (2.13) of \cite{hw:97}). Replacing $\sqrt{n} (S_n^{-1/2}- \textrm{I}_d )$ in the expression of $V_{n,1,1}(t)$ with the right-hand side of (\ref{snminushalb}), we obtain
\[
V_{n,1,1}(t) = - \frac{1}{2} \, \e^{-\|t\|^2/2} \left(\frac{1}{n} \sum_{j=1}^n \e^{t^\top X_j} X_j \right)^\top \frac{1}{\sqrt{n}} \sum_{j=1}^n \left(X_jX_j^\top - \textrm{I}_d\right) t + o_\PP(1).
\]
Replacing $n^{-1}\sum_{j=1}^n \exp(t^\top X_j) X_j$ with its expectation $\BE [\exp(t^\top X_1)X_1]$ $ = \exp(\|t\|^2/2) t$ means adding a term that is
$o_\PP(1)$ in $\LL$ which, upon combining with the non-negligible term of $V_{n,1,3}$, yields the assertion.
\bewend

\begin{prop} \label{asywn1}
We have
\[
W_{n,1}(t) =   \frac{1}{2} \, \frac{1}{\sqrt{n}} \sum_{j=1}^n t^\top \left( X_jX_j^\top - \textrm{I}_d \right) t + o_\PP(1),
\]
where $W_{n,1}$ is given in (\ref{wn1}) and $o_\PP(1)$ is with respect to $\LL$.
\end{prop}

\noindent {\sc Proof.} We have
$
W_{n,1}(t) = - W_{n,1,1}(t) + W_{n,1,2}(t) + W_{n,1,3}(t),
$
where
\begin{eqnarray*}
W_{n,1,1}(t) & = & \e^{\|t\|^2/2} \, \frac{1}{\sqrt{n}} \sum_{j=1}^n \sin\left(t^\top X_j \right) t^\top \left(S_n^{-1/2}- \textrm{I}_d \right) X_j,\\
W_{n,1,2}(t) & = & \e^{\|t\|^2/2} \, \frac{1}{\sqrt{n}} \sum_{j=1}^n \sin\left(t^\top X_j \right) t^\top \left(S_n^{-1/2}- \textrm{I}_d \right) \overline{X}_n,\\
W_{n,1,3}(t) & = & \e^{\|t\|^2/2} \, \frac{1}{\sqrt{n}} \sum_{j=1}^n \sin\left(t^\top X_j \right) t^\top \overline{X}_n.\\
\end{eqnarray*}
By complete analogy with the reasoning given in the proof of Proposition \ref{asyvn1}, we have $\|W_{n,1,2}\|^2_{\LL}  = o_\PP(1)$, and it is readily seen that also $\|W_{n,1,3}\|^2_{\LL} = o_\PP(1)$. Finally, by making use of (\ref{snminushalb}) and proceeding as in the proof of Proposition
\ref{asyvn1} (notice that $\BE[ \sin(t^\top X_1)X_1] = \exp(-\|t\|^2/2)t$), the assertion follows. \bewend



\section*{Acknowledgements}
M.D. Jim\'enez-Gamero was partially supported by
grant  MTM2014-55966-P of the Spanish Ministry of Economy and Competitiveness. Simos Meintanis was partially supported
by grant Nr.11699 of the Special Account for Research Grants (E${\rm{\Lambda}}$KE) of the National and
Kapodistrian University of Athens.

\end{document}